\documentclass[12pt,leqno]{amsart}

\usepackage{amssymb,amsmath,graphicx,amsfonts,euscript}
\usepackage{color}

\newtheorem{thm}{Theorem}[section]

\newtheorem{prop}[thm]{Proposition}
\newtheorem{define}[thm]{Definition}
\newtheorem{rem}[thm]{Remark}

\newtheorem{lemma}[thm]{Lemma}

\numberwithin{equation}{section}

\begin{document}
\title[Generalized surface quasi-geostrophic equations]
{Generalized surface quasi-geostrophic equations with singular velocities}
\author[Chae, Constantin, C\'{o}rdoba, Gancedo and Wu]{Dongho Chae$^{1}$, Peter Constantin$^{2}$, Diego C\'{o}rdoba$^{3}$, Francisco Gancedo$^{2}$, Jiahong Wu$^{4}$}

\address{$^1$ Department of Mathematics,
Sungkyunkwan University,
Suwon 440-746, Korea}

\address{$^2$Department of Mathematics,
 University of Chicago,
 5734 S. University Avenue,
Chicago, IL 60637, USA.}

\address{$^3$Instituto de Ciencias Matematicas,
Consejo Superior de Investigaciones Cientificas (CSIC),
28006 Madrid,
SPAIN}

\address{$^4$Department of Mathematics,
Oklahoma State University,
401 Mathematical Sciences,
Stillwater, OK 74078, USA.}

\email{chae@skku.edu}
\email{const@cs.uchicago.edu}
\email{dcg@icmat.es}
\email{fgancedo@math.uchicago.edu}
\email{jiahong@math.okstate.edu}

\begin{abstract}
This paper establishes several existence and uniqueness results for two families of active scalar equations with velocity fields determined by the scalars through very singular integrals. The first family is a generalized surface quasi-geostrophic (SQG) equation with the velocity field $u$ related to the scalar $\theta$ by $u=\nabla^\perp\Lambda^{\beta-2}\theta$, where $1<\beta\le 2$ and $\Lambda=(-\Delta)^{1/2}$ is the Zygmund operator. The borderline case $\beta=1$ corresponds to the SQG equation and the situation is more singular for $\beta>1$. We obtain the local existence and uniqueness of classical solutions, the global existence of weak solutions and the local existence of patch type solutions. The second family is a dissipative active scalar
equation with $u=\nabla^\perp (\log(I-\Delta))^\mu\theta$ for $\mu>0$, which is at least logarithmically more singular than the velocity in the first family. We prove that this family with any fractional dissipation possesses a unique local smooth solution for any given smooth data. This result for the second family constitutes a first step towards resolving the global regularity issue recently proposed by K. Ohkitani \cite{Oh}.
\end{abstract}

\maketitle

\section{Introduction}
\label{intr}
\setcounter{equation}{0}

This paper studies solutions of generalized surface quasi-geostrophic (SQG) equations
with velocity fields given by more singular integral operators than the Riesz transforms. Recall the inviscid SQG equation
\begin{align}
\begin{split}\label{sqg}
&\partial_t \theta + u\cdot\nabla \theta=0, \\
&u =\nabla^\perp \psi \equiv(-\partial_{x_2}, \partial_{x_1}) \psi, \quad \Lambda \psi =\theta,
\end{split}
\end{align}
where $\Lambda=(-\Delta)^{1/2}$ is the Zygmund operator, $\theta=\theta(x,t)$ is a scalar function, $u$ denotes the 2D velocity field and $\psi$ the stream function. Clearly, $u$ can be represented in terms of the Riesz transforms of $\theta$, namely
$$
u = (-\mathcal{R}_2, \mathcal{R}_1)\theta \equiv (-\partial_{x_2}\Lambda^{-1}, \partial_{x_1}\Lambda^{-1}) \theta.
$$
(\ref{sqg}), its counterpart with fractional dissipation and several closely related generalizations have recently been investigated very extensively and significant progress has been made on fundamental issues concerning solutions of these equations (see, e.g. \cite{AbHm}-\cite{Cham}, \cite{CMZ1}-\cite{KNV}, \cite{Li}-\cite{Sta}, \cite{WaJi}-\cite{Zhou2}).

\vskip .1in
Our goal here is to understand solutions of the SQG type equations with velocity fields determined by even more singular integral operators. Attention is focused on two generalized SQG equations. The first one assumes the form
\begin{align}
\begin{split}\label{active}
& \partial_t \theta + u\cdot\nabla \theta =0, \\
& u= \nabla^\perp \psi, \quad \Delta \psi =\Lambda^\beta\theta,
\end{split}
\end{align}
where $\beta$ is a real parameter satisfying  $1<\beta \le 2$. Here the spatial domain is either the whole plane $\mathbb{R}^2$ or the 2D periodic box $\mathbb{T}^2$ and the fractional Laplacian operator $(-\Delta)^\alpha$ is defined through the Fourier transform
$$
\widehat{(-\Delta)^\alpha f }(\xi) =|\xi|^{2 \alpha} \widehat{f}(\xi).
$$
The borderline case $\beta=1$ of (\ref{active}) is the SQG equation (\ref{sqg}), while (\ref{active}) with $\beta=0$ is the well-known 2D Euler vorticity equation with $\theta$ representing the vorticity (see e.g. \cite{MaBe}). The second generalized SQG equation under study is the dissipative active scalar equation
\begin{align}
\begin{split}\label{log2}
& \partial_t \theta + u\cdot\nabla \theta + \kappa (-\Delta)^\alpha\theta=0, \\
& u= \nabla^\perp \psi, \quad \psi =(\log(I-\Delta))^\mu\theta,
\end{split}
\end{align}
where $\kappa>0$, $\alpha>0$ and $\mu > 0$ are real parameters,
and $\left(\log(I-\Delta)\right)^\mu$ denotes the Fourier multiplier operator
defined by
$$
\widehat{\left(\log(I-\Delta)\right)^\mu f} (\xi) =  \left(\log(1+|\xi|^2)\right)^\mu \widehat{f} (\xi).
$$
(\ref{log2}) is closely related to (\ref{active}). In fact , both (\ref{active}) with $\beta=2$ and (\ref{log2}) with $\kappa=0$ and $\mu=0$ formally reduce to the trivial linear equation
$$
\partial_t \theta + \nabla^\perp \theta\cdot\nabla \theta =0 \quad\mbox{or}\quad \partial_t \theta =0.
$$
For $\mu>0$, the velocity field $u$ in (\ref{log2}) is at least logarithmically more singular than those in (\ref{active}).

\vskip .1in
We establish four main results for the existence and uniqueness of solutions to the equations defined in (\ref{active}) and in (\ref{log2}) with a given initial data
$$
\theta(x, 0) = \theta_0(x).
$$
We now preview these results. Our first main result establishes the local existence and uniqueness of smooth solutions to (\ref{active}) associated with any given smooth initial data.  More precisely, we have the following theorem.

\begin{thm} \label{local smooth}
Consider (\ref{active}) with $1< \beta\le 2$. Assume that $\theta_0\in H^m(\mathbb{R}^2)$ with $m \ge 4$. Then there exists $T = T(\|\theta_0\|_{H^m})>0$ such that (\ref{active})  has a unique solution $\theta$ on $[0,T]$. In addition, $\theta\in C([0,T]; H^m(\mathbb{R}^2))$.
\end{thm}

\begin{rem}
As mentioned previously, when $\beta=2$, $\psi =\theta$ and $u=\nabla^\perp \theta$ and (\ref{active}) reduces to the trivial equation
$$
\partial_t \theta =0 \quad \mbox{or}\quad \theta(x,t) =\theta_0(x).
$$
Therefore, (\ref{active}) with $\beta=2$ has a global steady-state solution.
\end{rem}

\vskip .1in
For $1<\beta <2$, the velocity $u$ is determined by a very singular integral of $\theta$ and $\nabla u$ is not known to be bounded in $L^\infty$. As a consequence, the nonlinear term can not be directly bounded. To deal with this difficulty, we rewrite the nonlinear term in the form of a commutator to explore the extra cancellation. In order to prove Theorem \ref{local smooth}, we need to derive a suitable commutator estimate (see Proposition \ref{comest} in Section \ref{local}).

\vskip .1in
Our second main result proves the local existence and uniqueness of smooth solutions to (\ref{log2}). In fact, the following theorem holds.

\begin{thm} \label{logthm}
Consider the active scalar equation (\ref{log2}) with $\kappa>0$, $\alpha>0$ and $\mu> 0$. Assume the initial data $\theta_0\in H^4(\mathbb{R}^2)$. Then there exists $T>0$ such that (\ref{log2}) has a unique solution $\theta \in C([0,T];H^4(\mathbb{R}^2))$.
\end{thm}

\vskip .1in
We remark that the velocity field $u$ in (\ref{log2}) is determined by
$$
u=\nabla^\perp (\log(I-\Delta))^\mu\theta \quad\mbox{with}\quad \mu>0
$$
which is even logarithmically more singular than that in (\ref{active}) with $\beta=2$, namely the trivial steady-state case. In a recent lecture \cite{Oh}, K. Ohkitani argued that (\ref{log2}) with $\kappa=0$ may be globally well-posed based on numerical computations. Theorem \ref{logthm} is a first step towards positively confirming his prediction.

\vskip .1in
Again the difficulty arises from the nonlinear term. In order to obtain a local (in time) bound for $\|\theta\|_{H^4}$, we need to rewrite the most singular part in the nonlinear term as a commutator. This commutator involves the logarithm of Laplacian and it appears that no $L^2$-bound for such commutator is currently available. By applying Besov space techniques, we are able to prove the following bound for such commutators.
\begin{prop} \label{comlog00}
Let $\mu\ge 0$. Let $\partial_x$ denote a partial derivative, either $\partial_{x_1}$ or $\partial_{x_2}$.  Then, for any $\delta>0$ and $\epsilon>0$, 
$$
\|\left[(\ln(I-\Delta))^\mu\partial_x, g\right] f \|_{L^2} \le C_{\mu,\epsilon,\delta} \left(1+ \left(\ln\left(1+ \frac{\|f\|_{\dot{H}^\delta}}{\|f\|_{L^2}}\right)\right)^{\mu}\right)\, \|f\|_{L^2}\, \|g\|_{H^{2+3\epsilon}},
$$
where $C_{\mu,\epsilon,\delta}$ is a constant depending on $\mu$, $\epsilon$ and $\delta$ only, $\dot{H}^\delta$ denotes the standard homogeneous Sobolev space and the brackets denote the commutator, namely
$$
\left[(\ln(I-\Delta))^\mu\partial_x, g\right] f = (\ln(I-\Delta))^\mu\partial_x(f g) - \left((\ln(I-\Delta))^\mu\partial_x f\right) \,g.
$$
\end{prop}

\vskip .1in
Our third main result assesses the global existence of weak solutions to (\ref{active}). Our consideration is restricted to the setting of periodic boundary conditions. The weak solution is essentially in the distributional sense and its precise definition is as follows.  $\mathbb{T}^2$ in the definition denotes the 2D periodic box.

\vskip .2in
\begin{define} \label{weak}
Let $T>0$. A function $\theta \in L^\infty([0,T];L^2(\mathbb{T}^2))$ is a weak solution of (\ref{active}) if, for any test function $\phi\in C^\infty_c([0,T)\times\mathbb{T}^2)$, the following integral equation holds,
\begin{equation}\label{weakeq}
\int_0^T \int_{\mathbb{T}^2} \theta\, (\partial_t\phi + u\cdot\nabla \phi)\, dx\,dt =
\int_{\mathbb{T}^2}\theta_0(x)\, \phi(x,0)\,dx.
\end{equation}
\end{define}

\vskip .1in
Although the velocity $u$ is more singular than the scalar $\theta$ and the nonlinear term above could not make sense, it is well defined due to a commutator hidden in the equation (see Section \ref{weee}). We prove that any mean-zero $L^2$ data leads to a global (in time) weak solution. That is, we have the following theorem.

\begin{thm} \label{weakthm}
Assume that $\theta_0\in L^2(\mathbb{T}^2)$ has mean zero, namely
$$
\int_{\mathbb{T}^2} \theta_0(x)\, dx =0.
$$
Then (\ref{active}) has a global weak solution in the sense of Definition \ref{weak}.
\end{thm}

This result is an extension of Resnick's work \cite{Res} on the inviscid SQG equation (\ref{sqg}).  However, for $1<\beta< 2$, the velocity is more singular and we need to write the nonlinear term as a commutator in terms of the stream function $\psi$. More details can be found in the proof of Theorem \ref{weakthm} in Section \ref{weee}.

\vskip .1in
Our last main result establishes the local well-posedness of the patch problem associated with the active scalar equation (\ref{active}). This result extends Gancedo's previous work for (\ref{active}) with $0<\beta\leq 1$ \cite{Gan}. Since $\beta$ is now in the range $(1,2)$, $u$ is given by a more singular integral and demands regular function and more sophisticated manipulation.  The initial data is given by
\begin{equation} \label{patch_data}
\theta_0(x) = \left\{
\begin{array}{ll}
\theta_1, \quad & x \in \Omega;\\
\theta_2, \quad & x\in \mathbb{R}^2\setminus \Omega,
\end{array}
\right.
\end{equation}
where $\Omega\subset \mathbb{R}^2$ is a bounded domain. We parameterize the boundary of $\Omega$ by $x=x_0(\gamma)$ with $\gamma\in \mathbb{T}=[-\pi, \pi]$ so that
$$
|\partial_\gamma x_0(\gamma)|^2 =A_0,
$$
where  $2\pi \sqrt{A_0}$ is the length of the contour. In addition, we assume that the curve $x_0(\gamma)$ does not cross itself and there is a lower bound on $|\partial_\gamma x_0(\gamma)|$, namely
\begin{equation} \label{sect}
\frac{|x_0(\gamma)-x_0(\gamma-\eta)|}{|\eta|} >0, \qquad \forall \gamma, \eta\in \mathbb{T}.
\end{equation}
Alternatively, if we define
\begin{equation} \label{fxg}
F(x)(\gamma,\eta,t) =\left\{
\begin{array}{l}
\displaystyle \frac{|\eta|}{|x(\gamma,t)-x(\gamma-\eta,t)|}, \quad \mbox{if}\, \eta\not = 0,\\ \\
\displaystyle \frac{1}{|\partial_\gamma x(\gamma, t)|}, \quad \mbox{if}\,  \eta = 0,
\end{array}
\right.
\end{equation}
then (\ref{sect}) is equivalent to
\begin{equation} \label{fin}
F(x_0)(\gamma,\eta,0)<\infty \quad \forall \gamma, \eta\in \mathbb{T}.
\end{equation}

\vskip .1in
The solution of (\ref{active}) corresponding to the initial data in (\ref{patch_data}) can be determined by studying the evolution of the boundary of the patch. As derived in \cite{Gan}, the parameterization $x(\gamma, t)$ of the boundary $\partial \Omega(t)$ satisfies
\begin{equation}\label{conto1}
\partial_t x(\gamma, t) = C_\beta (\theta_1-\theta_2) \int_T \frac{\partial_\gamma x(\gamma,t)- \partial_\gamma x(\gamma-\eta,t)}{|x(\gamma,t)-x(\gamma-\eta,t)|^\beta}\, d\eta,
\end{equation}
where $C_\beta$ is a constant depending on $\beta$ only. For $\beta\in (1,2)$, the integral on the right of (\ref{conto1}) is singular. Since the velocity in the tangential direction does not change the shape of the curve, we can modify (\ref{conto1}) in the tangential direction so that we get an extra cancellation. More precisely, we consider the modified equation
\begin{equation}\label{cde}
\partial_t x(\gamma, t) = C_\beta (\theta_1-\theta_2) \int_T \frac{\partial_\gamma x(\gamma,t)- \partial_\gamma x(\gamma-\eta,t)}{|x(\gamma,t)-x(\gamma-\eta,t)|^\beta}\, d\eta + \lambda(\gamma, t) \partial_\gamma x(\gamma,t)
\end{equation}
with $\lambda(\gamma,t)$ so chosen that
$$
\partial_\gamma x(\gamma,t) \cdot \partial_\gamma^2 x(\gamma,t) =0 \quad\mbox{or}\quad |\partial_\gamma x(\gamma,t)|^2 =A(t),
$$
where $A(t)$ denotes a function of $t$ only. A similar calculation as in \cite{Gan} leads to the following explicit formula for $\lambda(\gamma,t)$,
\begin{eqnarray}
\lambda(\gamma,t) &=& C\,\frac{\gamma+\pi}{2\pi} \int_\mathbb{T} \frac{\partial_\gamma x(\gamma,t)}{|\partial_\gamma x(\gamma,t)|^2}\cdot \partial_\gamma\left(\int_\mathbb{T} \frac{\partial_\gamma x(\gamma,t)- \partial_\gamma x(\gamma-\eta,t)}{|x(\gamma,t)-x(\gamma-\eta,t)|^\beta}\, d\eta\right)\,d\gamma \label{lameq}\\
&& - \,C\,\int_{-\pi}^\gamma \frac{\partial_\gamma x(\eta,t)}{|\partial_\gamma x(\eta,t)|^2}\cdot \partial_\eta\left(\int_\mathbb{T} \frac{\partial_\gamma x(\eta,t)- \partial_\gamma
x(\eta-\xi,t)}{|x(\eta,t)-x(\eta-\xi,t)|^\beta}\, d\xi\right)\,d\eta, \nonumber
\end{eqnarray}
where $C= C_\beta (\theta_1-\theta_2)$.

\vskip .1in
We establish the local well-posedness of the contour dynamics equation (CDE) given by (\ref{cde}) and (\ref{lameq}) corresponding to an initial contour
\begin{equation*}
x(\gamma,0) = x_0(\gamma)
\end{equation*}
satisfying (\ref{fin}). More precisely, we have the following theorem.
\begin{thm} \label{patch_local}
Let $x_0(\gamma)\in H^k(\mathbb{T})$ for $k\ge 4$ and $F(x_0)(\gamma,\eta,0)<\infty$ for any $\gamma, \eta\in \mathbb{T}$. Then there exists $T>0$ such that the CDE given by  (\ref{cde}) and (\ref{lameq}) has a solution $x(\gamma,t) \in C([0,T];H^k(\mathbb{T}))$ with $x(\gamma,0) = x_0(\gamma)$.
\end{thm}

\vskip .1in
This theorem is proven by obtaining an inequality of the form
$$
\frac{d}{dt}\left(\|x\|_{H^4} + \|F(x)\|_{L^\infty}\right) \le C \, \left(\|x\|_{H^4} + \|F(x)\|_{L^\infty}\right)^{9+\beta}.
$$
The ingredients involved in the proof include appropriate combination and cancellation of terms. The detailed proof is provided in Section \ref{patch}.

\vskip .4in
\section{Local smooth solutions}
\label{local}

This section proves Theorem \ref{local smooth}, which assesses the local (in time) existence and uniqueness of solutions to (\ref{active}) in $H^m$ with $m\ge 4$.

\vskip .1in
For $1<\beta\le 2$, the velocity $u$ is determined by a very singular integral of $\theta$ and the nonlinear term can not be directly bounded. To deal with this difficulty, we rewrite the nonlinear term in the form of a commutator to explore the extra cancellation. The following proposition provides a $L^2$-bound for the commutator
\begin{prop} \label{comest}
Let $s$ be a real number. Let $\partial_x$ denote a partial derivative, either $\partial_{x_1}$ or $\partial_{x_2}$. Then, 
$$
\|\left[\Lambda^s \partial_x, g\right] f\|_{L^2(\mathbb{R}^2)} \le
C\,\left(
\|\Lambda^s f\|_{L^2}\, \|\widehat{\Lambda g}(\eta)\|_{L^1} + C\,\|f\|_{L^2} \|\widehat{\Lambda^{1+s} g}(\eta)\|_{L^1} \right),
$$
where $C$ is a constant depending on $s$ only. In particular, by Sobolev embedding, 
for any $\epsilon>0$, there exists $C_\epsilon$ such that
$$
\|\left[\Lambda^s \partial_x, g\right] f\|_{L^2(\mathbb{R}^2)} \le
C_\epsilon\,\left(
\|\Lambda^s f\|_{L^2}\, \|g\|_{H^{2+\epsilon}} + \|f\|_{L^2}\, \|g\|_{H^{2+s+\epsilon}} \right).
$$
\end{prop}

\vskip .1in
Since this commutator estimate itself appears to be interesting, we provide a proof for this proposition.
\begin{proof}
The Fourier transform of $\left[\Lambda^s \partial_x, g\right] f$ is given by
\begin{equation} \label{ftran}
\widehat{\left[\Lambda^s \partial_x, g\right] f}(\xi) = \int_{\mathbb{R}^2} \left(
|\xi|^s \xi_j -|\xi-\eta|^s (\xi-\eta)_j\right) \widehat{f}(\xi-\eta)\, \widehat{g}(\eta)\, d\eta.
\end{equation}
where $j=1$ or $2$. It is easy to verify that, for any real number $s$,
\begin{equation} \label{basicin}
\left| |\xi|^s \xi_j -|\xi-\eta|^s (\xi-\eta)_j\right| \le C\, \max\{|\xi|^s, |\xi-\eta|^s\}\,|\eta|.
\end{equation}
In fact, we can write
\begin{eqnarray}
|\xi|^s \xi_j -|\xi-\eta|^s (\xi-\eta)_j &=& \int_0^1 \frac{d}{d\rho}\left(|A|^s A_j\right) \label{gdif} \\
&=& \int_0^1 (|A|^s \eta_j \,+ \, s|A|^{s-2} (A\cdot \eta) A_j) d\rho, \nonumber
\end{eqnarray}
where $A(\rho, \xi,\eta)= \rho \xi + (1-\rho)(\xi-\eta)$. Therefore,
$$
\left||\xi|^s \xi_j -|\xi-\eta|^s (\xi-\eta)_j \right| \le
(1+|s|) |\eta| \int_0^1 |A|^s \,d\rho.
$$
For $s\ge 0$, it is clear that
$$
|A|^s \le \max\{|\xi|^s, |\xi-\eta|^s\}.
$$
When $s<0$, $F(x)=|x|^s$ is convex and
$$
|A|^s = |\rho \xi + (1-\rho)(\xi-\eta)|^s \le \rho |\xi|^s + (1-\rho) |\xi-\eta|^s\le \max\{|\xi|^s, |\xi-\eta|^s\}.
$$
To obtain the bound in Proposition \ref{comest}, we first consider
the case when $s\ge 0$. Inserting (\ref{basicin}) in (\ref{ftran})
and using the basic inequality $|\xi|^s \le 2^{s-1} (|\xi-\eta|^s +
|\eta|^s)$, we have
\begin{eqnarray}
\left|\widehat{\left[\Lambda^s \partial_x, g\right] f}(\xi)\right| &\le&  C\,|\xi|^s \int_{\mathbb{R}^2} |\widehat{f}(\xi-\eta)|\, |\eta\widehat{g}(\eta)|\, d\eta \label{fff}\\
&& + \,C\,\int_{\mathbb{R}^2} ||\xi-\eta|^s\widehat{f}(\xi-\eta)|\, |\eta\widehat{g}(\eta)|\, d\eta \nonumber\\
&\le&  C\,\int_{\mathbb{R}^2} ||\xi-\eta|^s\widehat{f}(\xi-\eta)|\, |\eta\widehat{g}(\eta)|\, d\eta \nonumber\\
&& + \, C\,\int_{\mathbb{R}^2} |\widehat{f}(\xi-\eta)|\, ||\eta|^{1+s}\widehat{g}(\eta)|\, d\eta. \nonumber
\end{eqnarray}
By Plancherel's Theorem and Young's inequality for convolution,
$$
\|\left[\Lambda^s \partial_x, g\right] f\|_{L^2}
\le C\,\|\Lambda^s f\|_{L^2} \|\widehat{\Lambda g}(\eta)\|_{L^1} + C\,\|f\|_{L^2} \|\widehat{\Lambda^{1+s} g}(\eta)\|_{L^1}.
$$
Applying the embedding inequality
$$
\||\eta|^{1+s} \widehat{g}(\eta)\|_{L^1(\mathbb{R}^2)} \le C_\epsilon\,||g\|_{H^{2+s+\epsilon}(\mathbb{R}^2)},
$$
we have, for $s\ge 0$,
$$
\|\left[\Lambda^s \partial_x, g\right] f\|_{L^2(\mathbb{R}^2)} \le C_\epsilon\,\left(
\|\Lambda^s f\|_{L^2(\mathbb{R}^2)}\, \|g\|_{H^{2+\epsilon}(\mathbb{R}^2)}
+ \|f\|_{L^2(\mathbb{R}^2)}\, \|g\|_{H^{2+s+\epsilon}(\mathbb{R}^2)} \right).
$$
The case when $s<0$ is handled differently. We insert (\ref{gdif}) in (\ref{ftran}) and change
the order of integration to obtain
\begin{eqnarray*}
\widehat{\left[\Lambda^s \partial_x, g\right] f}(\xi) &=& H_1 + H_2,
\end{eqnarray*}
where
\begin{eqnarray}
H_1 &=&\int_0^1 \int_{\mathbb{R}^2} |A|^s  \widehat{f}(\xi-\eta)\, \eta_j\,\widehat{g}(\eta)\, d\eta\, d\rho, \label{H1} \\
H_2 &=& s \int_0^1 \int_{\mathbb{R}^2} |A|^{s-2} (A\cdot \eta) A_j \widehat{f}(\xi-\eta)\, \widehat{g}(\eta)
\,d\eta\, d\rho. \label{H2}
\end{eqnarray}
Using the fact that $F(x)=|x|^s$ with $s<0$ is convex, we have
\begin{eqnarray*}
|A|^s &=& |(\xi-\eta) + \rho \eta|^s = (1+\rho)^s \left|\frac{1}{1+\rho}(\xi-\eta) + \frac{\rho}{1+\rho} \eta\right|^s\\
&\le& (1+\rho)^s\left(\frac{1}{1+\rho}|\xi-\eta|^s + \frac{\rho}{1+\rho} |\eta|^s\right)\\
&=& (1+\rho)^{s-1} |\xi-\eta|^s + \rho (1+\rho)^{s-1} |\eta|^s.
\end{eqnarray*}
Inserting this inequality in (\ref{H1}), we obtain
\begin{eqnarray*}
|H_1| &\le& \int_0^1 (1+\rho)^{s-1} \,d\rho\, \int_{\mathbb{R}^2} ||\xi-\eta|^s\widehat{f}(\xi-\eta)|
\, |\eta\widehat{g}(\eta)|\, d\eta \\
&& + \int_0^1 \rho (1+\rho)^{s-1} \,d\rho\, \int_{\mathbb{R}^2} |\widehat{f}(\xi-\eta)| \,
||\eta|^{1+s} \,|\widehat{g}(\eta)|
\, |\eta.
\end{eqnarray*}
Applying Young's inequality for convolution, Plancherel's theorem and Sobolev's inequality, we have
\begin{eqnarray*}
\|H_1\|_{L^2} &\le& C\,\|\Lambda^s f\|_{L^2} \|\widehat{\Lambda g}(\eta)\|_{L^1} + C\,\|f\|_{L^2} \|\widehat{\Lambda^{1+s} g}(\eta)\|_{L^1}\\
&\le& C_\epsilon\, \|\Lambda^s f\|_{L^2}\, \|g\|_{H^{2+\epsilon}} +
C_\epsilon\,  \|f\|_{L^2}\, \|g\|_{H^{2+s+\epsilon}}.
\end{eqnarray*}
To bound $H_2$, it suffices to notice that
$$
|H_2| \le |s| \int_0^1 \int_{\mathbb{R}^2} |A|^s  |\widehat{f}(\xi-\eta)|\, |\eta \widehat{g}(\eta)|\, d\eta\, d\rho
$$
Therefore, $\|H_2\|_{L^2}$ admits the same bound as $\|H_1\|_{L^2}$. This completes
the proof of Proposition \ref{comest}.
\end{proof}

With this commutator estimate at our disposal, we are ready to prove Theorem \ref{local smooth}.

\vskip .1in
\begin{proof}[Proof of Theorem \ref{local smooth}] This proof provides a local (in time) {\it a priori} bound for $\|\theta\|_{H^m}$. Once the local bound is established, the construction of a local solution can be obtained through standard procedure such as the successive approximation. We shall omit the construction part to avoid redundancy.

\vskip .1in
We consider the case when $m=4$. The general case can be dealt with in a similar manner. By $\nabla \cdot u =0$,
$$
\frac12 \frac{d}{dt} \|\theta(\cdot,t)\|_{L^2}^2 = 0 \quad\mbox{or}\quad \|\theta(\cdot,t)\|_{L^2} = \|\theta_0\|_{L^2}.
$$
Let $\sigma$ be a multi-index with $|\sigma|=4$. Then,
$$
\frac12 \frac{d}{dt} \|D^\sigma \theta\|_{L^2}^2 = - \int D^\sigma \theta \, D^\sigma(u\cdot\nabla \theta)\,dx,
$$
where $\int$ means the integral over $\mathbb{R}^2$ and we shall
omit $dx$ when there is no confusion. Clearly, the right-hand side
can be decomposed into $I_1 +I_2 +I_3+I_4+ I_5$ with
\begin{eqnarray*}
I_1 &=& - \int D^\sigma \theta \,D^\sigma u \cdot\nabla \theta\,dx, \\
I_2 &=& - \sum_{|\sigma_1|=3, \sigma_1+\sigma_2=\sigma}\int D^\sigma \theta \, D^{\sigma_1}u\, \cdot D^{\sigma_2} \nabla \theta\,dx,\\
I_3 &=& - \sum_{|\sigma_1|=2, \sigma_1+\sigma_2=\sigma} \int D^\sigma \theta \, D^{\sigma_1}u\,\cdot D^{\sigma_2} \nabla \theta\,dx, \\
I_4&=& - \sum_{|\sigma_1|=1, \sigma_1+\sigma_2=\sigma}\int D^\sigma \theta \, D^{\sigma_1}u\, \cdot D^{\sigma_2} \nabla \theta\,dx,\\
I_5 &=& \int D^\sigma \theta \,u\, \cdot \nabla D^{\sigma}\theta\,dx.
\end{eqnarray*}
The divergence-free condition $\nabla \cdot u =0 $ yields $I_5=0$. We now estimate $I_1$. For $1<\beta<2$,
 $D^\sigma u=\nabla^\perp \Lambda^{-2+\beta} D^\sigma\theta$ with $|\sigma|=4$ can not bounded directly
 in terms of $\|\theta\|_{H^4}$. We rewrite $I_1$ as a commutator.
 For this we observe that for any skew-adjoint operator $A$ in $L^2$(i.e.
 $(Af,g)_{L^2}=-(f,Ag)_{L^2}$ for all $f,g\in L^2$) we have $\int f A(f) g \,dx= -\int f A(gf)\,dx $,
 and therefore
  \begin{equation}\label{commu}
  \int f A(f) g \,dx= -\frac12 \int \{ f A(gf)-fg A(f)\}\,dx= -\frac12 \int f[A,g]f\, dx.
 \end{equation}
Applying this fact to $I_1$ with $A:= \Lambda^{-2+\beta} \nabla
^\bot$, $f:=D^\sigma \theta$ and $g:=\nabla \theta$, one obtains
$$
I_1 = \frac12 \, \int D^\sigma \theta \,\left[\Lambda^{-2+\beta}
\nabla^\perp\cdot ,  \nabla\theta\right] D^\sigma \theta \,dx.
$$
By H\"{o}lder's inequality and Proposition \ref{comest} with $s=-2+\beta<0$, we have
$$
|I_1| \le C_\epsilon\, \|D^\sigma\theta\|_{L^2}\, \left(\|D^\sigma\theta\|_{L^2} +\| \Lambda^{-2+\beta}D^\sigma\theta\|_{L^2}\right) \|\theta\|_{H^{3+\epsilon}}\le C\,\|D^\sigma\theta\|_{L^2}\, \|\theta\|^2_{H^4}.
$$
The estimate for $I_2$ is easy. By H\"{o}lder's and Sobolev's inequalities,
$$
|I_2| \le C \|D^\sigma\theta\|_{L^2}\, \|\theta\|_{H^{2+\beta}}\, \|\theta\|_{H^4}.
$$
By H\'{o}lder's inequality and the Gagliardo-Nirenberg inequality,
\begin{eqnarray*}
|I_3| &\le&  C   \sum_{|\sigma_1|=2, \sigma_1+\sigma_2=4} \|D^\sigma\theta\|_{L^2} \, \|D^{\sigma_1} u\|_{L^4} \, \|D^{\sigma_2}\nabla \theta\|_{L^4} \\
&\le& C\, \|D^\sigma\theta\|_{L^2} \, \|\theta\|_{H^{\beta+1}}^{1/2}\, \|\theta\|_{H^{\beta+2}}^{1/2}\|\theta\|_{H^3}^{1/2}\, \|\theta\|_{H^4}^{1/2}\\
&\le& C\, \|D^\sigma\theta\|_{L^2} \, \|\theta\|_{H^3}\, \|\theta\|_{H^4}.
\end{eqnarray*}
By H\"{o}lder's and Sobolev's inequalities,
\begin{eqnarray*}
|I_4| &\le&  C \sum_{|\sigma_1|=1, \sigma_1+\sigma_2=4} \|D^\sigma\theta\|_{L^2} \,
\|D^{\sigma_1}u\|_{L^\infty}\, \|D^{\sigma_2} \nabla \theta\|_{L^2}\\
&\le&  C \|D^\sigma\theta\|_{L^2} \, \|\theta\|_{H^{\beta+2}} \, \|\theta\|_{H^4}.
\end{eqnarray*}
For $1<\beta< 2$, the bounds above yields
$$
\frac{d}{dt} \|\theta\|_{H^4}^2 \le C \, \|\theta\|_{H^4}^3.
$$
This inequality allows us to obtain a local (in time) bound for $\|\theta\|_{H^4}$.

\vskip .1in
In order to get uniqueness, one could check the evolution of two solution with the same initial data. With a similar approach, we find
$$
\frac{d}{dt}\|\theta_2-\theta_1\|_{H^1}\leq C(\|\theta_2\|_{H^4}+\|\theta_1\|_{H^4})\|\theta_2-\theta_1\|_{H^1}.
$$
An easy application of the Gronwall inequality provides $\theta_2=\theta_1$. This concludes the proof of Theorem \ref{local smooth}.
\end{proof}

\vskip .4in
\section{The case that is logarithmically beyond $\beta=2$}
\label{log}

This section focuses on the dissipative active scalar equation defined in
(\ref{log2}) and the goal is to prove Theorem \ref{logthm}.

\vskip .1in
As mentioned in the introduction, the major difficulty in proving this theorem is due to the fact that the velocity $u$ is determined by a very singular integral of $\theta$. To overcome this difficulty, we rewrite the nonlinear term in the form of a commutator to explore the extra cancellation. The commutator involves the logarithm of the Laplacian and we need a suitable bound for this type of commutator. The bound is stated in Proposition \ref{comlog00}, but we restated here.

\vskip .1in
\begin{prop} \label{comlog}
Let $\mu\geq 0$. Let $\partial_x$ denote a first partial, i.e., either $\partial_{x_1}$ or $\partial_{x_2}$.  Then, for any $\delta>0$ and $\epsilon>0$,
$$
\|\left[(\ln(I-\Delta))^\mu\partial_x, g\right] f \|_{L^2} \le C_{\mu,\epsilon,\delta} \left(1+ \left(\ln\left(1+ \frac{\|f\|_{\dot{H}^\delta}}{\|f\|_{L^2}}\right)\right)^{\mu}\right)\, \|f\|_{L^2}\, \|g\|_{H^{2+3\epsilon}},
$$
where $C_{\mu,\epsilon,\delta}$ is a constant depending on $\mu$, $\epsilon$ and $\delta$ only and $\dot{H}^\delta$ denotes the standard homogeneous Sobolev space.
\end{prop}

\begin{rem}
The constant $C_{\mu,\epsilon,\delta}$ approaches $\infty$ as $\delta\to 0$ or $\epsilon\to 0$. When $\mu=0$, the constant depends on $\epsilon$ only. 
\end{rem}

\vskip .1in
We shall also make use of the following lemma that bounds the $L^2$-norm of the logarithm of function.
\begin{lemma} \label{logl2}
Let $\mu\ge 0$ be a real number. Then, for any $\delta>0$, 
\begin{equation} \label{logl2in}
\|\left(\ln(I-\Delta)\right)^\mu f\|_{L^2} \le C_{\mu,\delta} \|f\|_{L^2} \left(\ln\left(1+ \frac{\|f\|_{\dot{H}^\delta}}{\|f\|_{L^2}}\right)\right)^\mu.
\end{equation}
where $C_{\mu,\delta}$ is a constant depending on $\mu$ and $\delta$ only.
\end{lemma}

\vskip .1in
In the rest of this section, we first prove Theorem \ref{logthm}, then Proposition \ref{comlog} and finally Lemma \ref{logl2}.

\vskip .1in
\begin{proof}[Proof of Theorem \ref{logthm}]
The proof obtains a local {\it a priori} bound for $\|\theta\|_{H^4}$. Once the local bound is at our disposal, a standard approach such as the successive approximation can be employed to provide a complete proof for the local existence and uniqueness. Since this portion involves no essential difficulties, the details will be omitted.

\vskip .1in
To establish the local $H^4$ bound, we start with the $L^2$-bound. By $\nabla \cdot u=0$,
$$
\frac12 \frac{d}{dt} \|\theta\|^2_{L^2} + \kappa \|\Lambda^\alpha\theta\|_{L^2}^2 =0 \quad\mbox{or}\quad \|\theta(\cdot, t)\|_{L^2} \le \|\theta_0\|_{L^2}.
$$
Now let $\sigma$ be a multi-index with $|\sigma|=4$. Then,
\begin{eqnarray}
\frac12 \frac{d}{dt} \|D^\sigma \theta\|_{L^2}^2 + \kappa \|\Lambda^\alpha D^\sigma\theta\|_{L^2}^2 &=& - \int D^\sigma\theta \, D^\sigma(u\cdot\nabla \theta) \,dx \label{bright}\\
 &=& J_1+J_2+J_3+J_4 + J_5, \nonumber
\end{eqnarray}
where
\begin{eqnarray*}
J_1 &=& - \int D^\sigma \theta \,D^\sigma u \cdot\nabla \theta\,dx, \\
J_2 &=& - \sum_{|\sigma_1|=3, \sigma_1+\sigma_2=\sigma}\int D^\sigma \theta \, D^{\sigma_1}u\,\cdot D^{\sigma_2} \nabla \theta\,dx,\\
J_3 &=& - \sum_{|\sigma_1|=2, \sigma_1+\sigma_2=\sigma} \int D^\sigma \theta \, D^{\sigma_1}u\,\cdot D^{\sigma_2} \nabla \theta\,dx, \\
J_4&=& - \sum_{|\sigma_1|=1, \sigma_1+\sigma_2=\sigma}\int D^\sigma \theta \, D^{\sigma_1}u\,\cdot D^{\sigma_2} \nabla \theta\,dx,\\
J_5 &=& \int D^\sigma \theta \,u\, \cdot \nabla D^{\sigma}\theta\,dx.
\end{eqnarray*}
By $\nabla \cdot u =0$, $J_5=0$. To bound $J_1$, we write it as a
commutator integral.  Applying (\ref{commu}) with $A:=\nabla ^\bot
(\log(I-\Delta))^\mu $, $f:=D^\sigma \theta $ and $g:= \nabla
\theta$, we have
\begin{eqnarray*}
J_1=\frac12 \int  D^\sigma \theta
\left[\left(\log(I-\Delta)\right)^\mu\nabla^\perp \cdot , \nabla
\theta\right]  D^\sigma \theta\, dx.
\end{eqnarray*}
By H\"{o}lder's inequality and Proposition \ref{comlog},
\begin{eqnarray*}
|J_1| &\le&  C\, \|D^\sigma \theta\|_{L^2} \|\left[\left(\log(I-\Delta)\right)^\mu\nabla^\perp\cdot, \nabla \theta\right] D^\sigma \theta\|_{L^2} \\
&\le& C\, \|D^\sigma \theta\|^2_{L^2}\, \|\nabla \theta\|_{H^{2+ \epsilon}} \left(1+ \left(\ln(1+ \|D^\sigma \theta\|_{H^\delta} )\right)^{\mu}\right)\\
&\le& C_\epsilon\, \|D^\sigma \theta\|^2_{L^2}\, \|\theta\|_{H^{3+\epsilon}}\, \left(\ln(1+ \|\theta\|_{H^{4+\delta}})\right)^\mu.
\end{eqnarray*}
Applying H\"{o}lder's inequality, Lemma \ref{logl2} and the Sobolev embedding
\begin{equation}\label{lin}
H^{1+\epsilon}(\mathbb{R}^2) \hookrightarrow L^\infty(\mathbb{R}^2), \quad \epsilon>0
\end{equation}
we obtain
\begin{eqnarray*}
|J_2| &\le& C \sum_{|\sigma_1|=3, \sigma_1+\sigma_2=4} \|D^\sigma \theta\|_{L^2} \, \|D^{\sigma_1}u\|_{L^2} \, \|D^{\sigma_2} \nabla \theta\|_{L^\infty} \\
&\le& C_\epsilon\, \|D^\sigma \theta\|^2_{L^2} \, \left(\ln(1+ \|\theta\|_{H^{4+\delta}})\right)^\mu\, \|\theta\|_{H^{3+\epsilon}}.
\end{eqnarray*}
To bound $J_3$, we first apply H\"{o}lder's inequality to obtain
$$
|J_3| \le C \sum_{|\sigma_1|=2, \sigma_1+\sigma_2=4} \|D^\sigma \theta\|_{L^2} \,\|D^{\sigma_1}u\|_{L^4} \, \|D^{\sigma_2} \nabla \theta\|_{L^4}.
$$
By the Sobolev inequality
$$
\|f\|_{L^4(\mathbb{R}^2)} \le C\, \|f\|^{1/2}_{L^2(\mathbb{R}^2)}\, \|\nabla f\|^{1/2}_{L^2(\mathbb{R}^2)}
$$
and applying Lemma \ref{logl2}, we have
\begin{eqnarray*}
|J_3| &\le& C \sum_{|\sigma_1|=2, \sigma_1+\sigma_2=4}\|D^\sigma \theta\|_{L^2}\, \|D^{\sigma_1}u\|^{1/2}_{L^2}\, \|\nabla D^{\sigma_1}u\|^{1/2}_{L^2}\, \|D^{\sigma_2} \nabla \theta\|^{1/2}_{L^2}\, \|\nabla D^{\sigma_2} \nabla \theta\|^{1/2}_{L^2}\\
&\le& C\, \|D^\sigma \theta\|_{L^2}\, \|\theta\|^2_{H^4} \, \left(\ln(1+ \|\theta\|_{H^{4+\delta}})\right)^\mu.
\end{eqnarray*}
By H\"{o}lder's inequality, (\ref{lin}) and Lemma \ref{logl2},
\begin{eqnarray*}
|J_4| &\le& C \sum_{|\sigma_1|=1, \sigma_1+\sigma_2=4} \|D^\sigma \theta\|_{L^2}\, \|D^{\sigma_1}u \|_{L^\infty}\, \|D^{\sigma_2} \nabla \theta\|_{L^2} \\
&\le& C \sum_{|\sigma_1|=1, \sigma_1+\sigma_2=4} \|D^\sigma \theta\|_{L^2}\, \|D^{\sigma_1}u \|_{H^{1+\epsilon}}\, \|D^{\sigma_2} \nabla \theta\|_{L^2} \\
&\le& C\, \|D^\sigma \theta\|_{L^2} \|\theta\|_{H^4}\, \|\theta\|_{H^{3+\epsilon}}\, \left(\ln(1+ \|\theta\|_{H^{3+\epsilon+\delta}})\right)^\mu.
\end{eqnarray*}
Let $0<\epsilon \le 1$ and $0<\delta <\alpha$. The estimates above on the right-hand side of (\ref{bright}) then implies that
$$
\frac12 \frac{d}{dt} \|D^\sigma \theta\|_{L^2}^2 + \kappa \|\Lambda^\alpha D^\sigma\theta\|_{L^2}^2 \le C\, \|\theta\|_{H^4}^3 \left(\ln(1+ \|\theta\|_{H^{4+\alpha}})\right)^\mu.
$$
This inequality is obtained for $|\sigma|=4$. Obviously, for $|\sigma|=1, 2$ and $3$,
the bound on the right remains valid. Therefore, if we sum the inequalities for $|\alpha|=1,2,3$ and $4$ and recalling , we have
$$
\frac12 \frac{d}{dt} \|\theta\|_{H^4}^2 + \kappa \|\theta\|^2_{H^{4+\alpha}} \le C\, \|\theta\|_{H^4}^3 \left(\ln(1+ \|\theta\|_{H^{4+\alpha}})\right)^\mu.
$$
The local (in time) {\it a priori} bound for $\|\theta\|_{H^4}$ then follows if we notice the simple inequality $(\ln(1+a))^\mu \le a$ for large $a>0$. This completes the proof of Theorem \ref{logthm}.
\end{proof}

\vskip .2in
We now present the proof of Proposition \ref{comlog}.

\begin{proof}[Proof of Proposition \ref{comlog}] The proof involves Besov spaces and related concepts such as the Fourier localization operator $\Delta_j$ for $j=-1,0,1,\cdots$ and the operator $S_j$.
These tools are now standard and can be found in several books, say \cite{Che}, \cite{Lem} and \cite{RuSi}. A self-contained quick introduction to the notation used in this proof can be found in \cite{ChCW}.

\vskip .1in
We start by identifying $L^2$ with the inhomogeneous Besov space $B^0_{2,2}$, namely
$$
\|f\|_{L^2}^2  = \sum_{j=-1}^\infty \|\Delta_j f\|_{L^2}^2.
$$
Let $N\ge 1$ be an integer to be determined later. We write
\begin{equation}\label{k1k2}
\|\left[(\ln(I-\Delta))^\mu\partial_x, g\right] f \|_{L^2}^2  = K_1 + K_2,
\end{equation}
where
\begin{eqnarray}
K_1 &=& \sum_{j=-1}^{N-1} \|\Delta_j \left[(\ln(I-\Delta))^\mu\partial_x, g\right] f\|_{L^2}^2, \label{k1}\\
K_2 &=& \sum_{j=N}^\infty \|\Delta_j \left[(\ln(I-\Delta))^\mu\partial_x, g\right] f\|_{L^2}^2. \label{k2}
\end{eqnarray}
Following Bony's notion of paraproducts,
$$
F G =\sum_{k} S_{k-1} F \, \Delta_k G + \sum_{k} \Delta_k F \, S_{k-1} G + \sum_{k} \Delta_k F  \widetilde{\Delta}_k G
$$
with $\widetilde{\Delta}_k =\Delta_{k-1} + \Delta_k + \Delta_{k+1}$, we have the decomposition
\begin{eqnarray}
\left[(\ln(I-\Delta))^\mu \partial_x, g\right] f &=& (\ln(I-\Delta))^\mu\partial_x(f g) -\left((\ln(I-\Delta))^\mu\partial_x f\right)\,g \label{comde}\\
&=& L_1 + L_2 + L_3, \nonumber
\end{eqnarray}
where
\begin{eqnarray*}
L_1 &=& \sum_{k} (\ln(I-\Delta))^\mu\partial_x \left(S_{k-1} f \, \Delta_k g \right) - S_{k-1} \left((\ln(I-\Delta))^\mu\partial_x f \right) \,\Delta_k g, \\
L_2 &=& \sum_{k} (\ln(I-\Delta))^\mu\partial_x \left(\Delta_k f \,  S_{k-1} g \right) -  \Delta_k\left((\ln(I-\Delta))^\mu\partial_x f \right) \,S_{k-1} g, \\
L_3 &=& \sum_{k} (\ln(I-\Delta))^\mu\partial_x \left(\Delta_k f \,  \widetilde{\Delta}_k g \right) -  \Delta_k\left((\ln(I-\Delta))^\mu\partial_x f \right) \,\widetilde{\Delta}_k g.
\end{eqnarray*}
Inserting the decomposition (\ref{comde}) in (\ref{k1}) and (\ref{k2}) yields the following corresponding decompositions in $K_1$ and $K_2$,
$$
K_1  \le K_{11} + K_{12} + K_{13}, \qquad K_2  \le K_{21} + K_{22} + K_{23}
$$
with
\begin{eqnarray*}
&& K_{11} = \sum_{j=-1}^{N-1} \|\Delta_j L_1\|_{L^2}^2, \quad K_{12} = \sum_{j=-1}^{N-1} \|\Delta_j L_2\|_{L^2}^2, \quad K_{13} = \sum_{j=-1}^{N-1} \|\Delta_j L_3\|_{L^2}^2,\\
&& K_{21} =  \sum_{j=N}^\infty \|\Delta_j L_1\|_{L^2}^2, \quad K_{22} = \sum_{j=N}^\infty  \|\Delta_j L_2\|_{L^2}^2, \quad K_{23} =  \sum_{j=N}^\infty \|\Delta_j L_3\|_{L^2}^2.
\end{eqnarray*}

\vskip .1in
Attention is now focused on bounding these terms and we start with $K_{11}$. When $\Delta_j$ is applied to $L_1$, the summation over $k$ in $L_1$ becomes a finite summation for $k$ satisfying $|k-j|\le 3$, namely
$$
\Delta_j L_1 = \sum_{|k-j|\le 3} \Delta_j \left((\ln(I-\Delta))^\mu\partial_x \left(S_{k-1} f \, \Delta_k g \right) - S_{k-1} \left((\ln(I-\Delta))^\mu\partial_x f \right) \,\Delta_k g\right).
$$
For the sake of brevity, we shall just estimate the representative term with $k=j$ in $\Delta_j L_1$. The treatment of the rest of the terms satisfying $|k-j|\le 3$ is similar and yields the same bound. Therefore,
$$
\|\Delta_j L_1\|_{L^2} \le C\,\|\Delta_j\left((\ln(I-\Delta))^\mu\partial_x \left(S_{j-1} f \, \Delta_j g \right) - S_{j-1} \left((\ln(I-\Delta))^\mu\partial_x f \right) \,\Delta_j g\right)\|_{L^2}.
$$
Without loss of generality, we set $\partial_x =\partial_{x_1}$. By Plancherel's theorem,
$$
\|\Delta_j L_1\|_{L^2}^2 \le C\,\left\|\Phi_j(\xi) \int_{\mathbb{R}^2} \left(H(\xi) - H(\xi-\eta)\right) \widehat{S_{j-1}f}(\xi-\eta)\,\widehat{\Delta_j g}(\eta)\, d\eta \right\|^2_{L^2},
$$
where $\Phi_j$ denotes the symbol of $\Delta_j$, namely $\widehat{\Delta_j f}(\xi) =\Phi_j(\xi) \widehat{f}(\xi)$ and
$$
H(\xi) = \left((\ln(1+|\xi|^2)\right)^\mu \xi_1.
$$
To further the estimate, we first invoke the inequality
\begin{eqnarray*}
&& |H(\xi)-H(\xi-\eta)| \\
&& \qquad \le |\eta|\,\left(\left(\ln(1+\max\{|\xi|^2, |\xi-\eta|^2\})\right)^\mu + \mu \left(\ln(1+\max\{|\xi|^2, |\xi-\eta|^2\})\right)^{\mu-1}\right).
\end{eqnarray*}
Clearly, the first term on the right-hand side dominates. We assume, without loss of generality, that
\begin{equation}\label{dein}
 |H(\xi)-H(\xi-\eta)| \le C\, |\eta|\,\left(\ln(1+\max\{|\xi|^2, |\xi-\eta|^2\})\right)^\mu.
\end{equation}
Noticing that
$$
\mbox{supp}\, \Phi_j,\,\ \mbox{supp}\,\widehat{\Delta_j g}\subset \{\xi \in \mathbb{R}^2: \, 2^{j-1} \le |\xi| < 2^{j+1}\},
$$
we have, for $-1 \le j\le N-1$,
\begin{eqnarray}
\|\Delta_j L_1\|_{L^2}^2 &\le& C\,\left\|\Phi_j(\xi) \int_{\mathbb{R}^2} \left(\ln(1+\max\{|\xi|^2, |\xi-\eta|^2\})\right)^\mu \right. \label{j111}\\
&& \qquad\qquad \qquad  \left. |\widehat{S_{j-1}f}(\xi-\eta)|\,|\eta \widehat{\Delta_j g}(\eta)|\, d\eta \right\|^2_{L^2} \nonumber \\
&\le& C\,\left(\ln(1+ 2^{2N})\right)^{2\mu} \left\|\Phi_j(\xi) \int_{\mathbb{R}^2} |\widehat{S_{j-1}f}(\xi-\eta)|\,|\eta \widehat{\Delta_j g}(\eta)|\, d\eta \right\|^2_{L^2}  \nonumber\\
&\le& C\,\left(\ln(1+ 2^{2N})\right)^{2\mu} \left\|\int_{\mathbb{R}^2} |\widehat{S_{j-1}f}(\xi-\eta)|\,|\eta \widehat{\Delta_j g}(\eta)|\, d\eta \right\|^2_{L^2}. \nonumber
\end{eqnarray}
By Young's inequality for convolution,
$$
\|\Delta_j L_1\|_{L^2}^2 \le  C\,\left(\ln(1+ 2^{2N})\right)^{2\mu} \|\widehat{S_{j-1}f}\|_{L^2}^2 \|\eta \widehat{\Delta_j g}(\eta)\|_{L^1}^2.
$$
By Plancherel's theorem and H\"{o}lder's inequality, for any $\epsilon>0$,
$$
\|\widehat{S_{j-1}f}\|_{L^2} = \|S_{j-1}f\|_{L^2} \le \|f\|_{L^2}, \quad \|\eta \widehat{\Delta_j g}(\eta)\|_{L^1} \le C_\epsilon \|\Lambda^{2+\epsilon} \Delta_j g\|_{L^2}
$$
Therefore,
\begin{eqnarray}
K_{11} &\le&  C_\epsilon\,\left(\ln(1+ 2^{2N})\right)^{2\mu} \|f\|_{L^2}^2 \sum_{j=-1}^{N-1} \|\Lambda^{2+\epsilon} \Delta_j g\|_{L^2} \nonumber\\
&\le&  C_\epsilon\,\left(\ln(1+ 2^{2N})\right)^{2\mu}\, \|f\|_{L^2}^2 \, \|g\|^2_{H^{2+\epsilon}}. \label{k11b}
\end{eqnarray}
We now estimate $K_{12}$. As in $\Delta_j L_1$, we have
$$
\Delta_j L_2 = \sum_{|k-j|\le 3} \Delta_j \left((\ln(I-\Delta))^\mu\partial_x \left(\Delta_k f \,  S_{k-1} g \right) -  \Delta_k\left((\ln(I-\Delta))^\mu\partial_x f \right) \,S_{k-1} g \right).
$$
It suffices to estimate the representative term with $k=j$. As in the estimate of $\Delta_j L_1$, we have
\begin{eqnarray*}
\|\Delta_j L_2\|_{L^2}^2 &\le& C\,\left(\ln(1+ 2^{2N})\right)^{2\mu} \left\|\int_{\mathbb{R}^2} |\widehat{\Delta_j f}(\xi-\eta)|\,|\eta \widehat{S_{j-1} g}(\eta)|\, d\eta \right\|^2_{L^2} \\
&\le& C\,\left(\ln(1+ 2^{2N})\right)^{2\mu} \|\widehat{\Delta_j f}\|_{L^2}^2 \|\eta \widehat{S_{j-1} g}(\eta)\|_{L^1}^2\\
&\le& C\, \left(\ln(1+ 2^{2N})\right)^{2\mu} \|\Delta_j f\|_{L^2}^2 \, \|g\|^2_{H^{2+\epsilon}}.
\end{eqnarray*}
Therefore,
\begin{eqnarray}
K_{12} &\le&  C\,\left(\ln(1+ 2^{2N})\right)^{2\mu}  \sum_{j=-1}^{N-1} \|\Delta_j f\|_{L^2}^2 \, \|g\|^2_{H^{2+\epsilon}} \nonumber \\
 &\le& C\,\left(\ln(1+ 2^{2N})\right)^{2\mu}  \|f\|_{L^2}^2 \|g\|^2_{H^{2+\epsilon}}. \label{k12b}
\end{eqnarray}
$K_{13}$ involves the interaction between high frequencies of $f$ and $g$ and the estimate is slightly more complicated. First we notice that
$$
\Delta_j L_3  = \sum_{k\ge j-1} \Delta_j \left((\ln(I-\Delta))^\mu\partial_x (\Delta_k f \,  \widetilde{\Delta}_k g) -  \Delta_k\left((\ln(I-\Delta))^\mu\partial_x f \right) \,\widetilde{\Delta}_k g \right).
$$
Applying Plancherel's theorem and invoking (\ref{dein}), we find
\begin{eqnarray}
\|\Delta_j L_3\|_{L^2}^2 &\le& \sum_{k\ge j-1} \Big\|\Delta_j \left((\ln(I-\Delta))^\mu\partial_x (\Delta_k f \,  \widetilde{\Delta}_k g)  \right.\\
&& \qquad\qquad\quad   -  \left. \Delta_k\left((\ln(I-\Delta))^\mu\partial_x f \right) \,\widetilde{\Delta}_k g \right)\Big\|_{L^2}^2 \label{jl31} \nonumber\\
&\le& C \sum_{k\ge j-1}  \Big\|\Phi_j(\xi) \int_{\mathbb{R}^2} \left(\ln(1+\max\{|\xi|^2, |\xi-\eta|^2\})\right)^\mu \nonumber\\
&&\hskip 1.5in   \times |\widehat{\Delta_k f}(\xi-\eta)|\,|\eta \widehat{\widetilde{\Delta}_k g}(\eta)|\, d\eta  \Big\|_{L^2}^2. \nonumber
\end{eqnarray}
Since $\Phi_j$ is supported on $\{\xi \in \mathbb{R}^2: \, 2^{j-1} \le |\xi| < 2^{j+1}\}$ and $\widehat{\Delta_k f}$ is on $\{\xi \in \mathbb{R}^2: \, 2^{k-1} \le |\xi| < 2^{k+1}\}$, we have, for $k\ge j-1$,
\begin{eqnarray*}
\left(\ln(1+\max\{|\xi|^2, |\xi-\eta|^2\})\right)^\mu &\le& \left(\ln(1+\max\{2^{2j+2}, 2^{2(k+1)}\}\right)^\mu\\
 &\le& \left(\ln(1+ 2^{2k+4})\right)^\mu.
\end{eqnarray*}
Therefore,
\begin{eqnarray*}
\|\Delta_j L_3\|_{L^2}^2 &\le& C\,
\sum_{k\ge j-1} \left(\ln(1+ 2^{2k+4})\right)^{2\mu}
\Big\|\Phi_j(\xi) \int_{\mathbb{R}^2} |\widehat{\Delta_k f}(\xi-\eta)|\,|\eta \widehat{\widetilde{\Delta}_k g}(\eta)|\, d\eta  \Big\|_{L^2}^2.
\end{eqnarray*}
When $\eta$ is in the support of $\widehat{\widetilde{\Delta}_k g}$, $|\eta|$ is comparable to $2^{k}$ and $|\eta|^{2\epsilon}\sim 2^{2\epsilon k}$. Using this fact and Young's inequality for convolution, we have
\begin{eqnarray*}
\|\Delta_j L_3\|_{L^2}^2 &\le& C\, \sum_{k\ge j-1} \left(\ln(1+ 2^{2k+4})\right)^{2\mu} 2^{-2\epsilon k}\,
\left\|\int_{\mathbb{R}^2} |\widehat{\Delta_k f}(\xi-\eta)|\,||\eta|^{1+2\epsilon} \widehat{\widetilde{\Delta}_k g}(\eta)|\, d\eta  \right\|_{L^2}^2 \\
&\le& C\, \sum_{k\ge j-1} \left(\ln(1+ 2^{2k+4})\right)^{2\mu} 2^{-2\epsilon k}\,
\|\widehat{\Delta_k f}\|_{L^2}^2 \,\||\eta|^{1+2\epsilon} \widehat{\widetilde{\Delta}_k g}(\eta)\|_{L^1}^2.
\end{eqnarray*}
Using the fact that
$$
\left(\ln(1+ 2^{2k+4})\right)^{2\mu} 2^{-\epsilon k} \le C_\epsilon, \qquad \||\eta|^{1+2\epsilon} \widehat{\widetilde{\Delta}_k g}(\eta)\|_{L^1} \le C_\epsilon \|g\|_{H^{2+3\epsilon}},
$$
we obtain
$$
\|\Delta_j L_3\|_{L^2}^2  \le C_\epsilon\, \|g\|^2_{H^{2+3\epsilon}} \,\sum_{k\ge j-1} 2^{-\epsilon k} \|\Delta_k f\|_{L^2}^2.
$$
Therefore,
\begin{eqnarray}
K_{13} &=& \sum_{j=-1}^{N-1} \|\Delta_j L_3\|_{L^2}^2 \nonumber\\
&\le&  C_\epsilon\, \|g\|^2_{H^{2+3\epsilon}} \,\sum_{j=-1}^{N-1} 2^{-\epsilon j} \sum_{k\ge j-1} 2^{-\epsilon (k-j)} \|\Delta_k f\|_{L^2}^2 \nonumber\\
&\le& C_\epsilon\, \|g\|^2_{H^{2+3\epsilon}} \, \|f\|_{L^2}^2. \label{k13b}
\end{eqnarray}
We now turn to $K_{21}$. $\Delta_j L_1$ is bounded differently. As in (\ref{j111}), we have
\begin{eqnarray*}
\|\Delta_j L_1\|_{L^2}^2  &\le& C\, \left\|\Phi_j(\xi) \int_{\mathbb{R}^2} \left(\ln(1+\max\{|\xi|^2, |\xi-\eta|^2\})\right)^\mu \right.\\
&& \qquad\qquad \qquad \left. |\widehat{S_{j-1}f}(\xi-\eta)|\,|\eta\widehat{\Delta_j g}(\eta)|\, d\eta \right\|^2_{L^2}.
\end{eqnarray*}
Since $\mbox{supp}\, \Phi_j,\,\ \mbox{supp}\,\widehat{\Delta_j g}\subset \{\xi \in \mathbb{R}^2: \, 2^{j-1} \le |\xi| < 2^{j+1}\}$,
we have
$$
\left(\ln(1+\max\{|\xi|^2, |\xi-\eta|^2\})\right)^\mu \le C\, \left(\ln(1+2^{2j})\right)^\mu
$$
and $\eta \in \mbox{supp}\,\widehat{\Delta_j g}$ indicates that $|\eta|$ is comparable with $2^{j}$. Therefore,
\begin{eqnarray*}
\|\Delta_j L_1\|_{L^2}^2  &\le& C\, \left(\ln(1+2^{2j})\right)^{2\mu} 2^{-2\epsilon j}\, \left\|\int_{\mathbb{R}^2}|\widehat{S_{j-1}f}(\xi-\eta)|\,||\eta|^{1+\epsilon} \widehat{\Delta_j g}(\eta)|\, d\eta \right\|^2_{L^2} \\
&\le& C\, \left(\ln(1+2^{2j})\right)^{2\mu} 2^{-2\epsilon j}\,\|\widehat{S_{j-1}f}\|_{L^2}^2 \||\eta|^{1+\epsilon} \widehat{\Delta_j g}(\eta)\|_{L^1}^2\\
&\le& C\, \left(\ln(1+2^{2j})\right)^{2\mu} 2^{-2\epsilon j}\,\|f\|_{L^2}^2 \|\Lambda^{2+2\epsilon} \Delta_j g\|_{L^2}^2.
\end{eqnarray*}
Therefore,
\begin{eqnarray}
K_{21} &=& \sum_{j=N}^\infty \|\Delta_j L_1\|_{L^2}^2 \nonumber \\
&\le& C\, \|f\|_{L^2}^2 \,  \sum_{j=N}^\infty \left(\ln(1+2^{2j})\right)^{2\mu} 2^{-2\epsilon j}\, \|\Lambda^{2+2\epsilon} \Delta_j g\|_{L^2}^2 \nonumber\\
&\le& C\,  \|f\|_{L^2}^2 \, \left(\ln(1+2^{2N})\right)^{2\mu} 2^{-2\epsilon N}\,\|g\|^2_{H^{2+2\epsilon}} \nonumber\\
&\le& C\,  \|f\|_{L^2}^2 \,\|g\|^2_{H^{2+2\epsilon}}. \label{k21b}
\end{eqnarray}
We now bound $K_{22}$. $\Delta_j L_2$ admits the following bound
\begin{eqnarray*}
\|\Delta_j L_2\|_{L^2}^2  &\le& C\, \left\|\Phi_j(\xi) \int_{\mathbb{R}^2} \left(\ln(1+\max\{|\xi|^2, |\xi-\eta|^2\})\right)^\mu \right.\\
&& \qquad\qquad \qquad \left. |\widehat{\Delta_j f}(\xi-\eta)|\,|\eta\widehat{S_{j-1} g}(\eta)|\, d\eta \right\|^2_{L^2}.
\end{eqnarray*}
Since  $\mbox{supp}\, \Phi_j\subset \{\xi \in \mathbb{R}^2: \, 2^{j-1} \le |\xi| < 2^{j+1}\}$ and $\mbox{supp}\,\widehat{S_{j-1} g}\subset \{\xi \in \mathbb{R}^2:  |\xi| < 2^{j}\}$, we still have
$$
\left(\ln(1+\max\{|\xi|^2, |\xi-\eta|^2\})\right)^\mu \le C\, \left(\ln(1+2^{2j})\right)^\mu.
$$
In contrast to the previous estimate on $\Delta_j L_1$, $\eta\in \widehat{S_{j-1} g}$ no longer implies that $|\eta|$ is comparable to $2^j$. However, any $\xi \in  \mbox{supp}\,\widehat{\Delta_j f}$ must have $|\xi|$ comparable to $2^j$. Therefore,
for any $\delta>0$,
\begin{eqnarray*}
\|\Delta_j L_2\|_{L^2}^2  &\le& C\,\left(\ln(1+2^{2j})\right)^{2\mu} 2^{-2\delta j}
\left\|\int_{\mathbb{R}^2} ||\xi-\eta|^{\delta}\widehat{\Delta_j f}(\xi-\eta)|\,|\eta\widehat{S_{j-1} g}(\eta)|\, d\eta \right\|^2_{L^2}\\
&\le& C\,\left(\ln(1+2^{2j})\right)^{2\mu} 2^{-2\delta j}  \||\xi-\eta|^{\delta}\widehat{\Delta_j f}(\xi-\eta)\|_{L^2}^2 \|\eta\widehat{S_{j-1} g}(\eta)\|_{L^1}^2 \\
&\le& C\,\left(\ln(1+2^{2j})\right)^{2\mu} 2^{-2\delta j} \|\Delta_j \Lambda^\delta f\|_{L^2}^2  \, \|g\|_{H^{2+\epsilon}}^2.
\end{eqnarray*}
Thus,
\begin{eqnarray}
K_{22} &\le& C\, \sum_{j=N}^\infty \left(\ln(1+2^{2j})\right)^{2\mu} 2^{-2\delta j} \|\Delta_j \Lambda^\delta f\|_{L^2}^2  \, \|g\|_{H^{2+\epsilon}}^2 \nonumber \\
&\le& C\, \left(\ln(1+2^{2N})\right)^{2\mu} 2^{-2\delta N}  \|g\|_{H^{2+\epsilon}}^2\sum_{j=N}^\infty \|\Delta_j \Lambda^\delta f\|_{L^2}^2 \nonumber\\
&\le& C\, \left(\ln(1+2^{2N})\right)^{2\mu} 2^{-2\delta N}  \|g\|_{H^{2+\epsilon}}^2
\|f\|^2_{H^\delta}.  \label{k22b}
\end{eqnarray}
The last term $K_{23}$ can be dealt with exactly as $K_{13}$. The bound for $K_{23}$ is \begin{equation} \label{k23b}
K_{23} \le C_\epsilon\, \|g\|^2_{H^{2+3\epsilon}} \, \|f\|_{L^2}^2.
\end{equation}

\vskip .1in
Collecting the estimates in (\ref{k11b}), (\ref{k12b}), (\ref{k13b}), (\ref{k21b}), (\ref{k22b}) and (\ref{k23b}), and inserting them in (\ref{k1k2}),  we obtain, for any integer $N>1$,
\begin{eqnarray*}
\|\left[(\ln(I-\Delta))^\mu\partial_x, g\right] f \|_{L^2}^2   &\le& C_\epsilon\,\left(\ln(1+ 2^{2N})\right)^{2\mu}\, \|f\|_{L^2}^2 \, \|g\|^2_{H^{2+\epsilon}}\\
&& + \, C_\epsilon\, \|f\|_{L^2}^2\, \|g\|^2_{H^{2+3\epsilon}}\\
&&  + \, C_\epsilon\, \left(\ln(1+2^{2N})\right)^{2\mu} 2^{-2\delta N}\,  \|f\|^2_{H^\delta}\, \|g\|_{H^{2+\epsilon}}^2.
\end{eqnarray*}
We now choose $N$ such that $2^{-2\delta N}\,  \|f\|^2_{H^\delta} \le C \|f\|_{L^2}^2$. In fact, we can choose 
\begin{equation}\label{Neq}
N = \left[\frac1{\delta}\log_2\frac{\|f\|_{H^\delta}}{\|f\|_{L^2}}\right].
\end{equation}
It then follows that
$$
\|\left[(\ln(I-\Delta))^\mu\partial_x, g\right] f \|_{L^2} \le C_{\mu,\epsilon,\delta} \left(1+ \left(\ln\left(1+ \frac{\|f\|_{H^\delta}}{\|f\|_{L^2}}\right)\right)^{\mu}\right)\, \|f\|_{L^2}\, \|g\|_{H^{2+3\epsilon}},
$$
where $C_{\mu,\epsilon,\delta}$ is a constant depending on $\mu$, $\epsilon$ and $\delta$ only. It is easy to see that the inhomogeneous Sobolev norm $\|f\|_{H^\delta}$ can be replaced by the homogeneous norm $\|f\|_{\dot{H}^\delta}$. This completes the proof of Proposition \ref{comlog}.
\end{proof}

\vskip .1in
Finally we prove Lemma \ref{logl2}.
\begin{proof}[Proof of Lemma \ref{logl2}]
Let $N \ge 1$ be an integer to be specified later. We write
$$
\|\left(\ln(I-\Delta)\right)^\mu f\|_{L^2}^2 = L_1 + L_2
$$
where
\begin{eqnarray*}
L_1 = \sum_{j=-1}^{N-1} \|\Delta_j \left(\ln(I-\Delta)\right)^\mu f\|_{L^2}^2,\quad
L_2 = \sum_{j=N}^{\infty} \|\Delta_j \left(\ln(I-\Delta)\right)^\mu f\|_{L^2}^2.
\end{eqnarray*}
According to Theorem 1.2 in \cite{ChCW}, we have, for $j\ge 0$,
$$
\|\Delta_j \left(\ln(I-\Delta)\right)^\mu f\|_{L^2} \le C\, \left(\ln(1 + 2^{2j})\right)^\mu
\|\Delta_j f\|_{L^2}.
$$
Clearly, for $j=-1$,
$$
\|\Delta_{-1} \left(\ln(I-\Delta)\right)^\mu f\|_{L^2} \le C\, \|\Delta_{-1} f\|_{L^2}.
$$
Therefore,
$$
L_1 \le C\, \left(\ln(1 + 2^{2N})\right)^{2\mu} \sum_{j=-1}^{N-1} \|\Delta_j f\|^2_{L^2} \le C\, \left(\ln(1 + 2^{2N})\right)^{2\mu} \|f\|_{L^2}^2.
$$
For any $\delta>0$,
\begin{eqnarray*}
L_2 &\le& \sum_{j=N}^{\infty} \left(\ln(1 + 2^{2j})\right)^{2\mu} 2^{-2\delta j} \,2^{2\delta j} \|\Delta_j f\|_{L^2}^2\\
&\le& \left(\ln(1 + 2^{2N})\right)^{2\mu} 2^{-2\delta N} \, \|f\|_{H^\delta}^2.
\end{eqnarray*}
Therefore,
$$
\|\left(\ln(I-\Delta)\right)^\mu f\|_{L^2}^2 \le C\, \left(\ln(1 + 2^{2N})\right)^{2\mu} \|f\|_{L^2}^2 + \left(\ln(1 + 2^{2N})\right)^{2\mu} 2^{-2\delta N} \, \|f\|_{H^\delta}^2.
$$
If we choose $N$ in a similar fashion as in (\ref{Neq}), we obtain the desired inequality (\ref{logl2in}). This completes the proof of Lemma \ref{logl2}.
\end{proof}

\vskip .4in
\section{Global weak solutions}
\label{weee}

This section establishes the global existence of weak solutions to (\ref{active}), namely Theorem \ref{weakthm}. The following commutator estimate will be used.

\begin{lemma}\label{cccm}
Let $s\ge 0$. Let $j=1$ or $2$. Then, for any $\epsilon>0$, there exists a constant $C$ depending on $s$ and $\epsilon$ such that
\begin{equation}\label{bin}
\|[\Lambda^s\partial_{x_j}, g] h \|_{L^2(\mathbb{T}^2)}\le C\,\left(\|h\|_{L^2}\,\|g\|_{H^{2+s+\epsilon}} +  \|\Lambda^s h\|_{L^2} \, \|g\|_{H^{2+\epsilon}}\right).
\end{equation}
\end{lemma}

\vskip .1in
Although the lemma is for the periodic setting, it can be proven in a similar manner as Proposition \ref{comest} and we thus omit its proof.

\vskip .1in
\begin{proof}[Proof of Theorem \ref{weakthm}]
The proof follows a standard approach, the Galerkin approximation.
Let $n>0$ be an integer and let $K_n$ denotes the subspace of $L^2(\mathbb{T}^2)$,
$$
K_n= \left\{ e^{im\cdot x}: \, m\not =0\,\,\mbox{and}\,\, |m|\le n \right\}.
$$
Let $\mathbb{P}_n$ be the projection onto $K_n$.
For each fixed $n$, we consider the solution of the projected equation,
\begin{eqnarray*}
&& \partial_t \theta_n + \mathbb{P}_n (u_n \cdot \nabla \theta_n) =0, \\
&& u_n = \nabla^\perp \Lambda^{-2+\beta} \theta_n,  \\
&& \theta_n(x,0) = \mathbb{P}_n \theta_0(x).
\end{eqnarray*}
This equation has a unique global solution $\theta_n$. Clearly, $\theta_n$ obeys the $L^2$ global bound
\begin{equation}\label{l2b}
\|\theta_n(\cdot,t)\|_{L^2} = \|\mathbb{P}_n \theta_0\|_{L^2} \le \|\theta_0\|_{L^2}.
\end{equation}
In addition, let $\psi_n$ be the corresponding stream function, namely $\Delta \psi_n =\Lambda^\beta \theta_n$. Then we have
\begin{eqnarray*}
\frac12 \frac{d}{dt} \left\|\Lambda^{1-\frac{\beta}{2}} \psi_n  \right\|_{L^2}^2
&=& - \int \psi_n \mathbb{P}_n (u_n \cdot \nabla \theta_n) \,dx \\
&=& - \int \psi_n \,u_n \cdot \nabla \theta_n \,dx.
\end{eqnarray*}
Noticing that $u_n = \nabla^\perp \psi_n$, we integrate by parts in the last term to obtain
$$
- \int \psi_n \,u_n \cdot \nabla \theta_n \,dx = \int \psi_n \,u_n \cdot \nabla \theta_n \,dx.
$$
Therefore,
\begin{equation}\label{phib}
\frac{d}{dt} \left\|\Lambda^{1-\frac{\beta}{2}} \psi_n  \right\|_{L^2}^2  =0 \quad\mbox{or}\quad  \left\|\Lambda^{1-\frac{\beta}{2}} \psi_n  \right\|_{L^2} \le \left\|\Lambda^{1-\frac{\beta}{2}} \psi_0  \right\|_{L^2}.
\end{equation}
Furthermore, for any $\phi \in H^{3+\epsilon}$ with $\epsilon>0$, we have
\begin{equation}\label{tt}
\int \partial_t \theta_n(x,t) \, \phi(x)\,dx = - \int (u_n \cdot \nabla \theta_n) \mathbb{P}_n \phi\,dx = \int \theta_n u_n \cdot \nabla \mathbb{P}_n \phi\,dx.
\end{equation}
On the one hand, $\theta_n = \Lambda^{2-\beta}\psi_n$ and
$$
\int \theta_n u_n \cdot \nabla \mathbb{P}_n \phi\,dx = \int \psi_n \, \Lambda^{2-\beta}\left(u_n \cdot \nabla \mathbb{P}_n \phi\right)\,dx = \int \psi_n \, \Lambda^{2-\beta}\left(\nabla^\perp\psi_n \cdot \nabla \mathbb{P}_n \phi\right)\,dx.
$$
On the other hand, $u_n =\nabla^\perp\psi_n$ and
$$
\int \theta_n u_n \cdot \nabla \mathbb{P}_n \phi\,dx = \int \theta_n \nabla^\perp\cdot( \psi_n \, \nabla \mathbb{P}_n \phi)\,dx = - \int \psi_n \nabla^\perp \Lambda^{2-\beta}\psi_n \cdot \nabla \mathbb{P}_n \phi\,dx.
$$
Thus,
$$
\int \theta_n u_n \cdot \nabla \mathbb{P}_n \phi\,dx =\frac12 \int \psi_n \left[\Lambda^{2-\beta}\nabla^\perp\cdot, \nabla \mathbb{P}_n \phi\right]\psi_n \,dx.
$$
It then follows from H\"{o}lder's inequality and Lemma \ref{cccm} that
\begin{eqnarray}
\left|\int \theta_n u_n \cdot \nabla \mathbb{P}_n \phi\,dx\right| &\le&  C\, \|\psi_n\|_{L^2}\, \|\psi_n\|_{H^{2-\beta}} \|\mathbb{P}_n \phi\|_{H^{3+\epsilon}} \label{triples}\\
&\le& C\,\|\Lambda^{-2+\beta}\theta_n\|_{L^2}\, \|\theta_n\|_{L^2} \,\|\phi\|_{H^{3+\epsilon}} \nonumber \\
&\le& C\, \|\theta_0\|_{L^2}^2 \, \|\phi\|_{H^{3+\epsilon}}  \nonumber
\end{eqnarray}
where we have used the fact that mean-zero functions in $L^2(\mathbb{T}^2)$ are also in $H^{-2+\beta}(\mathbb{T}^2)$. Therefore, by (\ref{tt}),
\begin{equation}\label{pttb}
\|\partial_t \theta_n \|_{H^{-3-\epsilon}} \le C\, \|\theta_0\|_{L^2}^2.
\end{equation}
The bounds in (\ref{l2b}), (\ref{phib}) and (\ref{pttb}), together with the compact embedding relation $L^2(\mathbb{T}^2)\hookrightarrow H^{-2+\beta}(\mathbb{T}^2)$ for $1<\beta<2$, imply that there exists $\theta\in C([0,T]; L^2(\mathbb{T}^2))$
such that
\begin{equation}\label{cov1}
\theta_n \rightharpoonup \theta \quad\mbox{in}\quad L^2, \qquad \psi_n \to \psi \quad\mbox{in}\quad L^2.
\end{equation}
In addition, because of the uniform boundedness of  $\|\theta_n\|_{L^2}$ and the embedding $L^2(\mathbb{T}^2)\hookrightarrow H^{-3-\epsilon}(\mathbb{T}^2)$, the Arzel\`{a}-Ascoli Theorem implies
\begin{equation}\label{cov2}
\lim_{n \to \infty} \sup_{t\in [0,T]} \left|\int (\theta_n(x,t)-\theta(x,t)) \phi(x) \,dx \right| \to 0,
\end{equation}
where $\phi\in H^{3+\epsilon}(\mathbb{T}^2)$.

\vskip .1in
The convergence in (\ref{cov1}) and (\ref{cov2}) allows us to prove that $\theta$ satisfies (\ref{weakeq}). Clearly, $\theta_n$ satisfies the integral equation
$$
\int_0^T \int_{\mathbb{T}^2} \theta_n\, (\partial_t\phi + u_n\cdot\nabla \mathbb{P}_n\phi)\, dx\,dt =\int_{\mathbb{T}^2} \mathbb{P}_n \theta_{0}(x)\, \phi(x,0)\,dx.
$$
It is easy to check that
\begin{eqnarray*}
&& \int_{\mathbb{T}^2} \mathbb{P}_n \theta_{0}(x)\, \phi(x,0)\,dx \to  \int_{\mathbb{T}^2} \theta_0(x)\, \phi(x,0)\,dx,
\end{eqnarray*}
and (\ref{cov2}) implies that, as $n\to \infty$,
\begin{eqnarray*}
&& \int_0^T \int_{\mathbb{T}^2} \theta_n\,\partial_t\phi \,dx dt \to \int_0^T \int_{\mathbb{T}^2} \theta\,\partial_t\phi \,dx dt.
\end{eqnarray*}
To show the convergence in the nonlinear term, we write
\begin{eqnarray*}
&& \int_0^T \int_{\mathbb{T}^2} \theta_n\,u_n\cdot \nabla\mathbb{P}_n\phi\, dx\,dt - \int_0^T \int_{\mathbb{T}^2} \theta\,u \cdot \nabla\phi\, dx\,dt \nonumber\\
&& \quad =\frac12 \int_0^T\int_{\mathbb{T}^2} \psi_n \left[\Lambda^{2-\beta}\nabla^\perp\cdot, \nabla \mathbb{P}_n \phi\right]\psi_n \,dx\,dt\nonumber\\
&& \quad\quad - \frac12 \int_0^T\int \psi \left[\Lambda^{2-\beta}\nabla^\perp\cdot, \nabla \phi\right]\psi \,dx\,dt \nonumber\\
&& \quad  = \frac12 \int_0^T\int \psi_n \left[\Lambda^{2-\beta}\nabla^\perp\cdot, \nabla (\mathbb{P}_n \phi-\phi)\right]\psi_n \,dx\,dt  \nonumber\\
&&\qquad + \frac12 \int_0^T \int_{\mathbb{T}^2}(\psi_n-\psi) \left[\Lambda^{2-\beta}\nabla^\perp\cdot, \nabla \phi\right]\psi_n \,dx\,dt \nonumber\\
&& \qquad + \frac12 \int_0^T \int_{\mathbb{T}^2}\psi \left[\Lambda^{2-\beta}\nabla^\perp\cdot, \nabla \phi\right](\psi_n-\psi) \,dx\,dt.
\label{thirdterm}
\end{eqnarray*}
In order to get the convergence for the first two terms above, we appeal to Lemma \ref{cccm} and the strong convergence of $\psi_n$ in $L^2$. Let us point out that in the last term for $\Lambda^{2-\beta}\psi_n$ we only have weak convergence in $L^2$ so we have to proceed in a different manner. We consider the following integral
\begin{eqnarray*}
Q_n(t)&=&\int_{\mathbb{T}^2}\psi \left[\Lambda^{2-\beta}\nabla^\perp\cdot, \nabla \phi\right](\psi_n-\psi) \,dx\\
&=& \sum_{k\neq 0}\widehat{\psi}(-k)(\left[\Lambda^{2-\beta}\nabla^\perp\cdot, \nabla \phi\right](\psi_n-\psi))\,\widehat{}\,\,(k),
\end{eqnarray*}
which is bounded by
$$
|Q_n(t)|\leq \left(\sum_{k\neq 0}\left||k|^{2-\beta}\widehat{\psi}(-k)\right|^2\right)^{1/2}
\left(\sum_{k\neq 0}\left||k|^{\beta-2}(\left[\Lambda^{2-\beta}\nabla^\perp\cdot, \nabla \phi\right](\psi_n-\psi))\,\widehat{}\,\,(k)\right|^2\right)^{1/2}.
$$
The first sum above is controlled by $\|\theta_0\|_{L^2}$. Using a similar notation as before, the coefficients in the second sum have the form
$$
|k|^{\beta-2}(\left[\Lambda^{2-\beta}\partial_x,\varphi\right](\psi_n-\psi))\,\widehat{}\,\,(k)
$$
where $\partial_x$ is either $\partial_{x_1}$ or $\partial_{x_2}$ and $\varphi$ is $\partial_{x}\phi$.
Since
\begin{align*}
(\left[\Lambda^{2-\beta}\partial_x,\varphi\right](\psi_n-\psi))\,\widehat{}\,\,(k)=\sum_{j}i(k_a|k|^{2-\beta}-(k-j)_a|k-j|^{2-\beta})
(\psi_n-\psi)\,\widehat{}\,\,(k-j)\widehat{\varphi}(j)
\end{align*}
for $a=1,2$, following the bounds in Section \ref{local} we obtain
\begin{align*}
\left|(\left[\Lambda^{2-\beta}\partial_x,\varphi\right](\psi_n-\psi))\,\widehat{}\,\,(k)\right|&\leq C\sum_{j}(|k|^{2-\beta}+|k-j|^{2-\beta})|(\psi_n-\psi)\,\widehat{}\,\,(k-j)||j||\widehat{\varphi}(j)|\\
&\leq C\sum_{j}(|k|^{2-\beta}+|j|^{2-\beta})|(\psi_n-\psi)\,\widehat{}\,\,(k-j)||j||\widehat{\varphi}(j)|.
\end{align*}
For $|k|\neq 0$, it yields
\begin{align*}
|k|^{\beta-2}\left|(\left[\Lambda^{2-\beta}\partial_x,\varphi\right](\psi_n-\psi))\,\widehat{}\,\,(k)\right|\leq C\sum_{j}|(\psi_n-\psi)\,\widehat{}\,\,(k-j)||j|(1+|j|^{2-\beta})|\widehat{\varphi}(j)|.
\end{align*}
The above bound provides
$$
|Q_n(t)| \leq C_{\epsilon}\|\theta_0\|_{L^2}\|\phi\|_{H^{5-\beta+\epsilon}}\|\psi_n-\psi\|_{L^2}
$$
for any $\epsilon>0$. It then follows from (\ref{cov1}) that $\lim_{n\to \infty} Q_n(t)=0$. The Dominated Convergence Theorem then leads to the desired convergence of the third term. Therefore, $\theta$ is a weak solution of \ref{active} in the sense of Definition \ref{weak}. This completes the proof of Theorem \ref{weakthm}.
\end{proof}

\vskip .4in
\section{Local existence for smooth patches}
\label{patch}

This section is devoted to proving Theorem \ref{patch_local}.

\begin{proof}[Proof of Theorem \ref{patch_local}] Since $\beta=2$ corresponds to the trivial steady-state solution, it suffices to consider the case when $1<\beta<2$. The major efforts are devoted to establishing {\it a priori} local (in time) bound for $\|x(\cdot,t)\|_{H^4} + \|F(x)\|_{L^\infty}(t)$ for $x$ satisfying the contour dynamics equation (\ref{cde}) and $F(x)(\gamma,\eta,t)$ defined in (\ref{fxg}).

\vskip .1in
This proof follows the ideas in Gancedo \cite{Gan}. The difference here is that the kernel in (\ref{cde}) is more singular but the function space concerned here is $H^4(\mathbb{T})$, which is more regular than in \cite{Gan} and compensates for the singularity of the kernel.

\vskip .1in
For notational convenience, we shall omit the coefficient $C_\beta (\theta_1-\theta_2)$ in the contour dynamics equation (\ref{cde}). In addition, the $t$-variable will sometimes be suppressed. We start with the $L^2$-norm. Dotting (\ref{cde}) by $x(\gamma,t)$ and integrating over $\mathbb{T}$, we have
$$
\frac12 \frac{d}{dt} \int_\mathbb{T} |x(\gamma,t)|^2 \,dx = I_1 + I_2,
$$
where
\begin{eqnarray*}
I_1 &=& \int_\mathbb{T} \int_\mathbb{T} x(\gamma,t)\cdot\frac{\partial_\gamma x(\gamma,t)- \partial_\gamma x(\gamma-\eta,t)}{|x(\gamma,t)-x(\gamma-\eta,t)|^\beta}\, d\eta\, d\gamma, \\
I_2 &=& \int_\mathbb{T} \lambda(\gamma)\,  x(\gamma,t)\cdot \partial_\gamma x(\gamma,t) \, d\gamma.
\end{eqnarray*}
$I_1$ is actually zero. In fact, by the symmetrizing process,
\begin{eqnarray*}
I_1 &=& \frac12 \int_\mathbb{T} \int_\mathbb{T} \frac{(x(\gamma)-x(\gamma-\eta))\cdot (\partial_\gamma x(\gamma)- \partial_\gamma x(\gamma-\eta))}{|x(\gamma)-x(\gamma-\eta)|^\beta}\, d\eta\,d\gamma \\
&=& \frac1{2(2-\beta)} \int_\mathbb{T} \int_\mathbb{T} \partial_\gamma\left( |x(\gamma)-x(\gamma-\eta)|^{2-\beta}\right)\,d\gamma d\eta\\
&=& 0.
\end{eqnarray*}
To bound $I_2$, we first apply H\"{o}lder's inequality to obtain
$$
|I_2| \le \|\lambda\|_{L^\infty}\, \|x\|_{L^2}\, \|\partial_\gamma x\|_{L^2}.
$$
By the representation of $\lambda$ in (\ref{lameq}) and using the fact that
$$
 \frac{1}{|\partial_\gamma x|^2} \le \|F(x)\|^2_{L^\infty}(t),
$$
we have
\begin{eqnarray*}
\|\lambda\|_{L^\infty} &\le& C\, \|F(x)\|^2_{L^\infty}(t) \int_\mathbb{T} |\partial_\gamma x|\left|\partial_\gamma \int_\mathbb{T} \frac{\partial_\gamma x(\gamma)- \partial_\gamma x(\gamma-\eta)}{|x(\gamma)-x(\gamma-\eta)|^\beta}\, d\eta\right|\,d\gamma \\
&=& C\, \|F(x)\|^2_{L^\infty}(t) \left(I_{21} + I_{22}\right),
\end{eqnarray*}
where
\begin{eqnarray*}
I_{21} &=& \int_\mathbb{T} |\partial_\gamma x|\int_\mathbb{T} \frac{|\partial^2_\gamma x(\gamma)- \partial_\gamma^2 x(\gamma-\eta)|}{|x(\gamma)-x(\gamma-\eta)|^\beta}\, d\eta\,d\gamma,\\
I_{22} &=& \int_\mathbb{T} |\partial_\gamma x|\int_\mathbb{T} \frac{|\partial_\gamma x(\gamma)- \partial_\gamma x(\gamma-\eta)|^2}{|x(\gamma)-x(\gamma-\eta)|^{\beta+1}}\, d\eta\,d\gamma.
\end{eqnarray*}
It is not hard to see that $I_{21}$ and $I_{22}$ can be bounded as follows.
\begin{eqnarray*}
I_{21} &\le& C \|F(x)\|^\beta_{L^\infty}(t)\, \|\partial_\gamma x\|_{L^2} \, \|\partial_\gamma^3 x\|_{L^2},\\
I_{22} &\le& C \|F(x)\|^{1+\beta}_{L^\infty}(t)\,\|\partial_\gamma^2 x\|_{L^2}^2\,\|\partial_\gamma x\|_{L^2}.
\end{eqnarray*}
Therefore,
$$
\frac{d}{dt} \|x\|_{L^2}^2 \le C\, \|F(x)\|^{3+\beta}_{L^\infty}(t)\, \|x\|_{H^3}^5.
$$
We now estimate $\|\partial_\gamma^4 x\|_{L^2}$.
$$
\frac12 \frac{d}{dt} \int_\mathbb{T} |\partial_\gamma^4 x|^2 \,d\gamma = I_3 + I_4,
$$
where
\begin{eqnarray*}
I_3 &=& C\, \int_\mathbb{T} \partial_\gamma^4 x(\gamma)\cdot  \partial_\gamma^4\int_\mathbb{T} \frac{(\partial_\gamma x(\gamma)- \partial_\gamma x(\gamma-\eta))}{|x(\gamma)-x(\gamma-\eta)|^\beta}\, d\eta\,d\gamma,\\
I_4 &=& \int_\mathbb{T} \partial_\gamma^4 x(\gamma)\cdot \partial_\gamma^4 (\lambda \partial_\gamma x)(\gamma)\, d\gamma.
\end{eqnarray*}
$I_3$ can be further decomposed into five terms, namely $I_3= I_{31} + I_{32} + I_{33} + I_{34} + I_{35}$, where
\begin{eqnarray*}
I_{31} &=& \int_\mathbb{T} \int_\mathbb{T} \partial_\gamma^4 x(\gamma)\cdot \frac{(\partial_\gamma^5 x(\gamma)- \partial_\gamma^5 x(\gamma-\eta))}{|x(\gamma)-x(\gamma-\eta)|^\beta}\, d\eta\,d\gamma, \\
I_{32} &=& 4 \int_\mathbb{T} \int_\mathbb{T}  \partial_\gamma^4 x(\gamma)\cdot  (\partial_\gamma^4 x(\gamma)- \partial_\gamma^4 x(\gamma-\eta)) \partial_\gamma(|x(\gamma)-x(\gamma-\eta)|^{-\beta})\, d\eta\,d\gamma, \\
I_{33} &=& 6 \int_\mathbb{T} \int_\mathbb{T}  \partial_\gamma^4 x(\gamma)\cdot  (\partial_\gamma^3 x(\gamma)- \partial_\gamma^3 x(\gamma-\eta)) \partial_\gamma^2(|x(\gamma)-x(\gamma-\eta)|^{-\beta})\, d\eta\,d\gamma, \\
I_{34} &=& 4 \int_\mathbb{T} \int_\mathbb{T}  \partial_\gamma^4 x(\gamma)\cdot  (\partial_\gamma^2 x(\gamma)- \partial_\gamma^2 x(\gamma-\eta)) \partial_\gamma^3(|x(\gamma)-x(\gamma-\eta)|^{-\beta})\, d\eta\,d\gamma, \\
I_{35} &=&  \int_\mathbb{T} \int_\mathbb{T}  \partial_\gamma^4 x(\gamma)\cdot  (\partial_\gamma x(\gamma)- \partial_\gamma x(\gamma-\eta)) \partial_\gamma^4(|x(\gamma)-x(\gamma-\eta)|^{-\beta})\, d\eta\,d\gamma.
\end{eqnarray*}
By symmetrizing, $I_{31}$ can be written as
\begin{eqnarray*}
I_{31} &=& \frac12 \int_\mathbb{T} \int_\mathbb{T} \frac{(\partial_\gamma^4 x(\gamma)- \partial_\gamma^4 x(\gamma-\eta))\cdot (\partial_\gamma^5 x(\gamma)- \partial_\gamma^5 x(\gamma-\eta))}{|x(\gamma)-x(\gamma-\eta)|^\beta}\, d\eta\,d\gamma \\
&=& \frac14 \int_\mathbb{T} \int_\mathbb{T} \frac{\partial_\gamma (|\partial_\gamma^4 x(\gamma)- \partial_\gamma^4 x(\gamma-\eta)|^2)}{|x(\gamma)-x(\gamma-\eta)|^\beta}\, d\eta\,d\gamma \\
&=& \frac{\beta}4 \int_\mathbb{T} \int_\mathbb{T} \frac{|\partial_\gamma^4 x(\gamma)- \partial_\gamma^4 x(\gamma-\eta)|^2(x(\gamma)-x(\gamma-\eta))\cdot(\partial_\gamma x(\gamma)-\partial_\gamma x(\gamma-\eta))}{|x(\gamma)-x(\gamma-\eta)|^{\beta+2}}\, d\eta\,d\gamma. \\
\end{eqnarray*}
Setting
$$
B(\gamma, \eta)= (x(\gamma)-x(\gamma-\eta))\cdot(\partial_\gamma x(\gamma)-\partial_\gamma x(\gamma-\eta))
$$
and using the fact that $\partial_\gamma x(\gamma) \cdot \partial_\gamma^2 x(\gamma)=0$, we have
$$
|I_{31}| \le C\, \|F(x)\|^{2+\beta}_{L^\infty}(t)\,\int_\mathbb{T} \int_\mathbb{T} |\partial_\gamma^4 x(\gamma)- \partial_\gamma^4 x(\gamma-\eta)|^2 \frac{B(\gamma, \eta) \eta^{-2} -\partial_\gamma x(\gamma) \cdot \partial_\gamma^2 x(\gamma)}{|\eta|^{\beta}}d\eta\,d\gamma.
$$
Using the bound that
$$
\left|B(\gamma, \eta) \eta^{-2} -\partial_\gamma x(\gamma) \cdot \partial_\gamma^2 x(\gamma)\right| \le C\, \|x\|_{C^3}^2 |\eta|,
$$
we obtain
$$
|I_{31}| \le C\, \|F(x)\|^{2+\beta}_{L^\infty}(t)\, \|x\|_{C^3}^2 \|x\|_{H^4}^2.
$$
To estimate of $I_{32}$, we realize that, after computing $\partial_\gamma(|x(\gamma)-x(\gamma-\eta)|^{-\beta})$, $I_{32}$ can be bounded in the same fashion as $I_{31}$. That is,
$$
|I_{32}| \le C\, \|F(x)\|^{2+\beta}_{L^\infty}(t)\, \|x\|_{H^4}^4.
$$
In order to estimate $I_{33}$, we further decompose it into three terms, $I_{33} = I_{331} + I_{332} + I_{333}$, where
\begin{eqnarray*}
I_{331} &=& C\, \int_\mathbb{T} \int_\mathbb{T}  \partial_\gamma^4 x(\gamma)\cdot  (\partial_\gamma^3 x(\gamma)- \partial_\gamma^3 x(\gamma-\eta))
\frac{D(\gamma, \eta)}{|x(\gamma)-x(\gamma-\eta)|^{2+\beta}}\, d\eta\,d\gamma, \\
I_{332} &=& C\, \int_\mathbb{T} \int_\mathbb{T}  \partial_\gamma^4 x(\gamma)\cdot  (\partial_\gamma^3 x(\gamma)- \partial_\gamma^3 x(\gamma-\eta))
\frac{|\partial_\gamma x(\gamma)- \partial_\gamma x(\gamma-\eta)|^2}{|x(\gamma)-x(\gamma-\eta)|^{2+\beta}}\, d\eta\,d\gamma, \\
I_{333} &=& C\, \int_\mathbb{T} \int_\mathbb{T}  \partial_\gamma^4 x(\gamma)\cdot  (\partial_\gamma^3 x(\gamma)- \partial_\gamma^3 x(\gamma-\eta))
\frac{B^2(\gamma, \eta)}{|x(\gamma)-x(\gamma-\eta)|^{4+\beta}}\, d\eta\,d\gamma\\
\end{eqnarray*}
with
$$
D(\gamma, \eta) = (x(\gamma)- x(\gamma-\eta))\cdot(\partial_\gamma^2 x(\gamma)- \partial_\gamma^2 x(\gamma-\eta)).
$$
It is not very difficult to see that
$$
|I_{331}|, |I_{332}|, |I_{333}| \le C\, \|F(x)\|^{2+\beta}_{L^\infty}(t)\, \|x\|_{H^4}^4.
$$
$I_{34}$ also admit similar bound. In $I_{35}$ one has to use identity
$$
\partial_\gamma x(\gamma)\cdot\partial_\gamma^4 x(\gamma)=3\partial_\gamma^2 x(\gamma)\cdot\partial_\gamma^3 x(\gamma)
$$ to find the same control. We shall not provide the detailed estimates since they can be obtained by modifying the lines in \cite{Gan}.  We also need to deal with $I_4$. To do so, we use the representation formula (\ref{lameq}) and obtain
$$
|I_4| \le C\|F(x)\|^{4+\beta}_{L^\infty}(t)\, \|x\|_{H^4}^5
$$
In summary, we have
\begin{equation} \label{xl4}
\frac{d}{dt} \|x\|_{H^4}^2 \le C\, \|F(x)\|^{4+\beta}_{L^\infty}(t)\, \|x\|_{H^4}^5.
\end{equation}
We now derive the estimate for $\|F(x)\|_{L^\infty}(t)$. For any $p>2$, we have
\begin{equation} \label{flp}
\frac{d}{dt} \|F(x)\|^p_{L^p}(t) \le p \int_\mathbb{T} \int_\mathbb{T} \left(\frac{|\eta|}{|x(\gamma)- x(\gamma-\eta)|}\right)^{p+1}\frac{|x_t(\gamma,t)- x_t(\gamma-\eta,t)|}{|\eta|}\,d\eta\,d\gamma.
\end{equation}
Invoking the contour dynamics equation (\ref{cde}), we have
\begin{eqnarray*}
x_t(\gamma)- x_t(\gamma-\eta) &=& I_5+I_6+I_7+I_8\\
&\equiv&\int_\mathbb{T}\left(\frac{\partial_\gamma x(\gamma)- \partial_\gamma x(\gamma-\xi)}{|x(\gamma)- x(\gamma-\xi)|^\beta}-\frac{\partial_\gamma x(\gamma)-\partial_\gamma x(\gamma-\xi)}{|x(\gamma-\eta)- x(\gamma-\eta-\xi)|^\beta}\right) d\xi\\
&& + \int_\mathbb{T}\frac{\partial_\gamma x(\gamma)- \partial_\gamma x(\gamma-\eta)+\partial_\gamma  x(\gamma-\eta-\xi)-\partial_\gamma x(\gamma-\xi)}{|x(\gamma-\eta)- x(\gamma-\eta-\xi)|^\beta} d\xi\\
&& + (\lambda(\gamma)-\lambda(\gamma-\eta))\partial_\gamma x(\gamma) + \lambda(\gamma-\eta)(\partial_\gamma x(\gamma)-\partial_\gamma x(\gamma-\eta)).
\end{eqnarray*}
Following the argument as in \cite{Gan}, we have
\begin{eqnarray*}
|I_5| &\le&  C\, \|F(x)\|^{2\beta}_{L^\infty}(t)\, \|x\|^{1+\beta}_{C^2}\,|\eta|,\\
|I_6| &\le&  C\, \|F(x)\|^{\beta}_{L^\infty}(t)\,\|x\|_{C^3}\,|\eta|,\\
|I_7| &\le&  C\, \|F(x)\|^{3+\beta}_{L^\infty}(t)\, \|x\|^4_{H^4}\,|\eta|,\\
|I_8| &\le&  C\, \|F(x)\|^{3+\beta}_{L^\infty}(t)\, \|x\|^4_{H^4}\,|\eta|.
\end{eqnarray*}
Inserting these estimates in (\ref{flp}), we find
$$
\frac{d}{dt} \|F(x)\|_{L^p}(t) \le C\, \|x\|_{H^4}^4 \, \|F(x)\|^{4+\beta}_{L^\infty}(t)\, \|F(x)\|_{L^p}(t).
$$
After integrating in time and taking the limit as $p\to \infty$, we obtain
$$
\frac{d}{dt} \|F(x)\|_{L^\infty}(t) \le C\, \|x\|_{H^4}^4 \, \|F(x)\|^{5+\beta}_{L^\infty}(t).
$$
Combining with (\ref{xl4}), we obtain
$$
\frac{d}{dt} \left(\|x\|_{H^4} +\|F(x)\|_{L^\infty}(t) \right) \le C\, \|x\|_{H^4}^4\,\|F(x)\|^{5+\beta}_{L^\infty}(t).
$$
This inequality would allow us to deduce a local (in time) bound for $\|x\|_{H^4}$. This completes the proof of Theorem \ref{patch_local}.
\end{proof}

\vskip .4in
\section*{Acknowledgements}
This work was partially completed when Chae, Constantin, Gancedo and Wu visited the Instituto de Ciencias Matem\'{a}ticas (ICMAT), Madrid, Spain in November, 2010 and they thank the ICMAT for support and hospitality. Chae's research was partially supported by NRF grant No.2006-0093854. Constantin's research was partially supported by NSF grant DMS 0804380. Cordoba and Gancedo were partially supported by the grant MTM2008-03754 of the MCINN (Spain) and the grant StG-203138CDSIF of the ERC. Gancedo was also partially supported by NSF grant DMS-0901810. Wu's research was partially supported by NSF grant DMS 0907913 and he thanks Professors Hongjie Dong, Susan Friedlander and
Vlad Vicol for discussions.


\vskip .4in
\begin{thebibliography}{99}
\bibitem{AbHm} H. Abidi and T. Hmidi, On the global well-posedness of the critical quasi-geostrophic equation, {\it SIAM J. Math. Anal. \bf 40} (2008), 167--185.

\bibitem{Bae} H. Bae, Global well-posedness of dissipative quasi-geostrophic equations in critical spaces. {\it Proc. Amer. Math. Soc. \bf 136} (2008),  257--261.

\bibitem{Ber} L. Berselli, Vanishing viscosity limit and long-time behavior for 2D quasi-geostrophic equations, {\it Indiana Univ. Math. J. \bf 51} (2002), 905-930.


\bibitem{Blu} W. Blumen, Uniform potential vorticity flow, Part I.
Theory of wave interactions and two-dimensional turbulence, {\it J.
Atmos. Sci.} {\bf 35} (1978), 774-783.

\bibitem{CaS} L. Caffarelli and L. Silvestre, An extension problem
related to the fractional Laplacian, {\it Comm. Partial Differential Equations \bf 32} (2007), 1245--1260.

\bibitem{CV} L. Caffarelli and A. Vasseur, Drift diffusion
equations with fractional diffusion and the quasi-geostrophic
equation, {\it Ann. of Math. \bf 171} (2010),  1903-1930.

\bibitem{CaFe} J. Carrillo and L. Ferreira, The asymptotic behaviour of subcritical dissipative quasi-geostrophic equations, {\it Nonlinearity \bf 21} (2008), 1001-1018.

\bibitem{CaCo} A. Castro and D. C\'ordoba, Infinite energy solutions of the
surface quasi-geostrophic equation, {\it Adv. Math. \bf 223} (2010), 120-173.

\bibitem{Ch} D. Chae, The quasi-geostrophic equation in the Triebel-Lizorkin spaces,
{\it Nonlinearity} {\bf 16} (2003), 479-495.

\bibitem{ChJDE} D. Chae, On the continuation principles for the Euler equations and the quasi-geostrophic equation, {\it J. Differential Equations \bf 227} (2006), 640--651.

\bibitem{Cha}D. Chae, On the regularity conditions for the
dissipative quasi-geostrophic equations, {\it SIAM J. Math. Anal.}
{\bf 37} (2006), 1649-1656.

\bibitem{Cha2}D. Chae, The geometric approaches to the possible singularities in the inviscid fluid flows, {\it J. Phys. A  \bf 41} (2008), 365501, 11 pp.

\bibitem{Cha4} D. Chae, On the a priori estimates for the Euler, the Navier-Stokes and the quasi-geostrophic equations, {\it Adv. Math. \bf 221} (2009), 1678--1702.

\bibitem{ChCW} D. Chae, P. Constantin and J. Wu, Inviscid models generalizing the 2D Euler and the surface quasi-geostrophic equations,  arXiv:1010.1506v1  [math.AP]  7 Oct 2010.

\bibitem{ChCW2} D. Chae, P. Constantin and J. Wu, Dissipative models generalizing the 2D Navier-Stokes and the surface quasi-geostrophic equations,  arXiv:1011.0171 [math.AP] 31 Oct 2010.

\bibitem{CCCF} D. Chae, A. C\'{o}rdoba, D. C\'{o}rdoba and M. Fontelos, Finite time singularities in a 1D model of the quasi-geostrophic equation, {\it Adv. Math. \bf 194} (2005),  203--223.

\bibitem{ChL} D. Chae and J. Lee, Global well-posedness in the super-critical
dissipative quasi-geostrophic equations, {\it Commun. Math. Phys.}
{\bf 233} (2003), 297-311.

\bibitem{Cham} D. Chamorro, Remarks on a fractional diffusion transport equation with applications to the critical dissipative quasi-geostrophic equation,   arXiv:1007.3919v2 [math.AP] 22 Oct 2010.
    
\bibitem{ChCz}C. Chan, M. Czubak and L. Silvestre, Eventual regularization of the slightly supercritical fractional Burgers equation, {\it Discrete Contin. Dyn. Syst. \bf 27} (2010),  847-861.

\bibitem{Che} J.-Y. Chemin, {\it Perfect Incompressible Fluids},
Oxford science publications, Oxford University Press, 1998.


\bibitem{CMZ1} Q. Chen, C. Miao and Z. Zhang,
A new Bernstein's inequality and the 2D dissipative
quasi-geostrophic equation, {\it Commun. Math. Phys. \bf 271}
(2007), 821-838.

\bibitem{Chen} Q. Chen and Z. Zhang,
Global well-posedness of the 2D critical dissipative
quasi-geostrophic equation in the Triebel-Lizorkin spaces, {\it
Nonlinear Anal. \bf 67} (2007), 1715-1725.

\bibitem{Con} P. Constantin, Euler equations, Navier-Stokes equations and turbulence.
Mathematical foundation of turbulent viscous flows,  1--43, Lecture
Notes in Math., 1871, Springer, Berlin, 2006.

\bibitem{CCW} P. Constantin, D. C\'{o}rdoba and J. Wu, On the critical dissipative
quasi-geostrophic equation, {\em Indiana Univ. Math. J.} {\bf 50}
(2001), 97-107.


\bibitem{CIW} P. Constantin, G. Iyer and J. Wu, Global regularity for a modified critical dissipative quasi-geostrophic equation, {\it Indiana Univ. Math. J.}
    {\bf 57} (2008),  2681-2692.

\bibitem{CLS}P. Constantin, M.-C. Lai, R. Sharma, Y.-H. Tseng and J. Wu, New numerical results for the surface quasi-geostrophic equation, submitted for publication.
\bibitem{CMT} P. Constantin, A. Majda, and  E. Tabak,
Formation of strong fronts in the 2-D quasi-geostrophic thermal
active scalar, {\it Nonlinearity} {\bf 7} (1994), 1495-1533.

\bibitem{CNS} P. Constantin, Q. Nie and N. Schorghofer,  Nonsingular
surface quasi-geostrophic flow, {\it Phys. Lett. A } {\bf  241}
(1998),  168--172.

\bibitem{CW5} P. Constantin and J. Wu, Behavior of solutions of 2D
quasi-geostrophic equations, {\it SIAM J. Math. Anal.} {\bf 30}
(1999), 937-948.

\bibitem{CWnew1} P. Constantin and J. Wu, Regularity of H\"{o}lder continuous
solutions of the supercritical quasi-geostrophic equation, {\it Ann.
Inst. H. Poincar\'{e} Anal. Non Lin\'{e}aire} {\bf 25} (2008), 1103-1110.

\bibitem{CWnew2} P. Constantin and J. Wu,  H\"{o}lder continuity of solutions
of supercritical dissipative hydrodynamic transport equation, {\it
Ann. Inst. H. Poincar\'{e} Anal. Non Lin\'{e}aire} {\bf 26} (2009), 159-180.

\bibitem{Cor} D. C\'{o}rdoba, Nonexistence of simple hyperbolic blow-up for the quasi-geostrophic equation, {\it Ann. of Math. \bf 148} (1998),  1135--1152.

\bibitem{CC} A. C\'{o}rdoba and D. C\'{o}rdoba, A maximum principle
applied to quasi-geostrophic equations, {\it Commun. Math. Phys.}
{\bf 249} (2004), 511-528.


\bibitem{CoFe1} D. C\'{o}rdoba and Ch. Fefferman, Behavior of several
two-dimensional fluid equations in singular scenarios, {\it Proc.
Natl. Acad. Sci. USA} {\bf  98}  (2001),  4311--4312.

\bibitem{CoFe2} D. C\'{o}rdoba and Ch. Fefferman, Scalars convected by
a two-dimensional incompressible flow, {\it Comm. Pure Appl. Math.}
{\bf 55} (2002),  255--260.

\bibitem{CoFe3} D. C\'{o}rdoba and Ch. Fefferman, Growth of solutions for QG and 2D
Euler equations, {\it J. Amer. Math. Soc.} {\bf  15} (2002),
665--670.

\bibitem{CFMR} D. C\'{o}rdoba, M. Fontelos, A. Mancho and J. Rodrigo, Evidence of singularities for a family of contour dynamics equations, {\it Proc. Natl. Acad. Sci. USA \bf 102} (2005),  5949--5952.

\bibitem{Dab} M. Dabkowski, Eventual regularity of the solutions to the
supercritical dissipative quasi-geostrophic equation, arXiv:1007.2970v1 [math.AP] 18 Jul 2010.

\bibitem{DHLY} J. Deng, T. Y. Hou, R. Li and X. Yu, Level set dynamics and the non-blowup of the 2D quasi-geostrophic equation, {\it Methods Appl. Anal. \bf 13} (2006),  157--180.


\bibitem{DoCh} B. Dong and Z. Chen, Asymptotic stability of the critical and super-critical dissipative quasi-geostrophic equation, {\it Nonlinearity \bf 19} (2006), 2919-2928.

\bibitem{Dong}H. Dong, Dissipative quasi-geostrophic equations in critical Sobolev spaces: smoothing effect and global well-posedness, {\it Discrete Contin. Dyn. Syst. \bf 26} (2010), 1197--1211.

\bibitem{DoDu} H. Dong and D. Du, Global well-posedness and a decay estimate for the critical dissipative quasi-geostrophic equation in the whole space, {\it Discrete Contin. Dyn. Syst. \bf 21} (2008),  1095--1101.

\bibitem{DoLi0} H. Dong and D. Li,  Finite time singularities for a class of generalized surface quasi-geostrophic equations, {\it Proc. Amer. Math. Soc. \bf 136}(2008),  2555--2563.

\bibitem{DoLi} H. Dong and D. Li, Spatial analyticity of the solutions to the subcritical dissipative quasi-geostrophic equations, {\it Arch. Ration. Mech. Anal. \bf 189} (2008),  131--158.
    
\bibitem{DoLi2} H. Dong and D. Li, On the 2D critical and supercritical dissipative quasi-geostrophic equation in Besov spaces, {\it J. Differential Equations \bf 248} (2010),  2684-2702.

\bibitem{DoPo} H. Dong and N. Pavlovic, A regularity criterion for the dissipation quasi-geostrophic equation,  {\it Ann. Inst. H. Poincar\'{e} Anal. Non Lin\'{e}aire} {\bf 26} (2009),  1607--1619.

\bibitem{DoPo2} H. Dong and N. Pavlovic,  Regularity criteria for the dissipative quasi-geostrophic equations in H\"{o}lder spaces, {\it Comm. Math. Phys. \bf 290} (2009),  801--812.

\bibitem{FPV} S. Friedlander, N. Pavlovic and V. Vicol, Nonlinear instability for the critically dissipative quasi-geostrophic equation, {\it Comm. Math. Phys. \bf 292} (2009), 797--810.

\bibitem{FrVi}S. Friedlander and V. Vicol, Global well-posedness for an advection-diffusion equation arising in magneto-geostrophic dynamics,     arXiv:1007.1211v1 [math.AP] 12 Jul 2010.

\bibitem{Gan}  F. Gancedo, Existence for the $\alpha$-patch model and the QG sharp front in Sobolev spaces, {\it Adv. Math. \bf 217} (2008),  2569-2598.

\bibitem{Gil} A.E. Gill, {\it Atmosphere-Ocean Dynamics}, Academic Press,
New York, 1982.

\bibitem{HaOg}H. Hayashi and T. Ogawa, $L^p-L^q$ type estimate for the fractional order Laplacian in the Hardy space and global existence of the dissipative quasi-geostrophic equation, {\it Adv. Differ. Equ. Control Process. \bf 5} (2010), 1-36.
    
\bibitem{HPGS} I. Held, R. Pierrehumbert, S. Garner, and K. Swanson,
Surface quasi-geostrophic dynamics, {\it J. Fluid Mech.} {\bf 282}
(1995), 1-20.



\bibitem{HmKe} T. Hmidi and S. Keraani, Global solutions of the super-critical 2D quasi-geostrophic equation in Besov spaces, {\it Adv. Math. \bf 214} (2007), 618--638.

\bibitem{HmKe2} T. Hmidi and S. Keraani, On the global well-posedness of the critical quasi-geostrophic equation, {\it SIAM J. Math. Anal. \bf 40} (2008), 167--185.
    
\bibitem{HoSh}T.Y. Hou and Z. Shi, Dynamic growth estimates of maximum vorticity for 3D incompressible Euler equations and the SQG model, 	arXiv:1011.5514v1 [math.AP] 24 Nov 2010.

\bibitem{Ju} N. Ju, The maximum principle and the global attractor for the dissipative
2D quasi-geostrophic equations,  {\it Commun. Math. Phys. \bf 255}
(2005), 161-181.

\bibitem{Ju2} N. Ju, Geometric constrains for global regularity of 2D quasi-geostrophic flows, {\it J. Differential Equations \bf 226} (2006), 54--79.

\bibitem{KhTi}B. Khouider and E. Titi, An inviscid regularization for the surface quasi-geostrophic equation, {\it Comm. Pure Appl. Math. \bf 61} (2008), 1331--1346.

\bibitem{Ki1} A. Kiselev, Some recent results on the critical surface quasi-geostrophic equation: a review, {\it Hyperbolic problems: theory, numerics and applications}, 105--122, Proc. Sympos. Appl. Math., 67, Part 1, AMS, Providence, RI, 2009.

\bibitem{Ki2} A. Kiselev, Regularity and blow up for active scalars, {\it Math. Model. Math. Phenom. \bf 5} (2010), 225--255.

\bibitem{Kinew1} A. Kiselev, Nonlocal maximum principles for active scalars, arXiv: 1009.0542v1 [math.AP] 2 Sep 2010.


\bibitem{KN1}A. Kiselev and F. Nazarov, Global regularity for the critical dispersive dissipative surface quasi-geostrophic equation, {\it Nonlinearity \bf 23} (2010), 549--554.

\bibitem{KN2}A. Kiselev and F. Nazarov, A variation on a theme of Caffarelli and Vasseur, {\it Zap. Nauchn. Sem. POMI \bf 370} (2010), 58--72.

\bibitem{KNV}A. Kiselev, F. Nazarov and A. Volberg, Global well-posedness for the
critical 2D dissipative quasi-geostrophic equation, {\it Invent.
Math. \bf 167} (2007), 445-453.


\bibitem{Lem}P.-G. Lemari\'{e}-Rieusset, ``Recent developments in the Navier-Stokes
problem", Chapman \& Hall/CRC, Boca Raton, FL, 2002.

\bibitem{Li} D. Li, Existence theorems for the 2D quasi-geostrophic equation with plane wave initial conditions, {\it Nonlinearity \bf 22} (2009), 1639--1651.

\bibitem{LiRo} D. Li and J. Rodrigo, Blow up for the generalized surface quasi-geostrophic equation with supercritical dissipation, {\it Comm. Math. Phys. \bf 286} (2009),  111--124.

\bibitem{LiZh}P. Li and Z.Zhai, Riesz transforms on Q-type spaces with application to quasi-geostrophic equation, 	arXiv:0907.0856v1 [math.AP] 5 Jul 2009.
    
\bibitem{Maj} A. Majda, {\it Introduction to PDEs and Waves for the
Atmosphere and Ocean}, Courant Lecture Notes {\bf 9}, Courant
Institute of Mathematical Sciences and American Mathematical
Society, 2003.

\bibitem{MaBe}A. Majda and A. Bertozzi,  {\it Vorticity and Incompressible Flow}, Cambridge University Press, 2002.

\bibitem{MaTa} A. Majda and E. Tabak, A two-dimensional model for quasigeostrophic flow: comparison with the two-dimensional Euler flow, {\it  Phys. D \bf 98} (1996), 515--522.

\bibitem{Man} A. Mancho, Numerical studies on the self-similar collapse of the $\alpha$-patches problem, 	arXiv:0902.0706v1 [math.AP] 4 Feb 2009.

\bibitem{Mar1}F. Marchand, Propagation of Sobolev regularity for the critical dissipative quasi-geostrophic equation, {\it Asymptot. Anal. \bf 49} (2006),  275--293.

\bibitem{Mar2}F. Marchand, Existence and regularity of weak solutions to the quasi-geostrophic equations in the spaces $L^p$ or $\dot H^{-1/2}$, {\it Comm. Math. Phys. \bf 277} (2008),  45--67.

\bibitem{Mar3}F. Marchand, Weak-strong uniqueness criteria for the critical quasi-geostrophic equation, {\it Phys. D  \bf 237} (2008),  1346--1351.

\bibitem{MarLr} F. Marchand and P.G. Lemari\'{e}-Rieusset, Solutions
auto-similaires non radiales pour l'\'{e}quation
quasi-g\'{e}ostrophique dissipative critique, {\it  C. R. Math.
Acad. Sci. Paris} {\bf  341}  (2005),  535--538.

\bibitem{May} R. May, Global well-posedness for a modified 2D dissipative quasi-geostrophic equation with initial data in the critical Sobolev space $H^1$,
    arXiv:0910.0998v1 [math.AP] 6 Oct 2009.

\bibitem{MayZ}R. May and E. Zahrouni, Global existence of solutions for subcritical quasi-geostrophic equations, {\it Commun. Pure Appl. Anal. \bf 7} (2008), 1179--1191.

\bibitem{MiXu} C. Miao and L. Xue,
Global wellposedness for a modified critical dissipative quasi-geostrophic equation,
arXiv:0901.1368v4 [math.AP] 18 Sep 2010.

\bibitem{MiXu2} C. Miao and L. Xue, On the regularity of a class of generalized quasi-geostrophic equations, arXiv:1011.6214v1 [math.AP] 29 Nov 2010.


\bibitem{Mi} H. Miura, Dissipative quasi-geostrophic equation for large initial data in the critical sobolev space, {\it Commun. Math. Phys. \bf 267} (2006), 141--157.

\bibitem{NiSc} C. Niche and M. Schonbek, Decay of weak solutions to the 2D dissipative quasi-geostrophic equation, {\it Comm. Math. Phys. \bf 276} (2007), 93--115.

\bibitem{Oh} K. Ohkitani, Dissipative and ideal surface
quasi-geostrophic equations, Lecture presented at ICMS, Edinburgh, September, 2010.

\bibitem{OhSa} K. Ohkitani and T. Sakajo, Oscillatory damping in long-time evolution of the surface quasi-geostrophic equations with generalised viscosity: a numerical study, preprint.

\bibitem{OhYa1} K. Ohkitani and M. Yamada, Inviscid and inviscid-limit behavior
of a surface quasigeostrophic flow, {\it Phys. Fluids } {\bf  9}
(1997), 876--882.


\bibitem{ReDr}J. Reinaud and D. Dritschel, Destructive interactions between two counter-rotating quasi-geostrophic vortices, {\it J. Fluid Mech. \bf 639} (2009), 195--211.

\bibitem{Res} S. Resnick, Dynamical problems in nonlinear advective
partial differential equations, Ph.D. thesis, University of Chicago,
1995.

\bibitem{Ro1} J. Rodrigo,  The vortex patch problem for the surface quasi-geostrophic equation, {\it Proc. Natl. Acad. Sci. USA  \bf 101} (2004),  2684--2686

\bibitem{Ro2} J. Rodrigo, On the evolution of sharp fronts for the quasi-geostrophic equation, {\it Comm. Pure Appl. Math. \bf 58} (2005),  821--866.

\bibitem{RuSi} T. Runst and W. Sickel, {\it Sobolev Spaces of Fractional Order, Nemytskij Operators, and Nonlinear Partial Differential Equations}, Walter de Gruyter \& Co., Berlin, 1996.

\bibitem{Sch} M. Schonbek and T. Schonbek, Asymptotic behavior to
dissipative quasi-geostrophic flows, {\it SIAM J. Math. Anal.} {\bf
35} (2003), 357-375.

\bibitem{Sch2} M. Schonbek and T. Schonbek, Moments and lower
bounds in the far-field of solutions to quasi-geostrophic flows,
{\it Discrete Contin. Dyn. Syst.} {\bf 13} (2005), 1277-1304.

\bibitem{Si} L. Silvestre, Eventual regularization for the slightly supercritical quasi-geostrophic equation, {\it Ann. Inst. H. Poincar\'{e} Anal. Non Lin\'{e}aire \bf 27} (2010), no. 2, 693--704.

\bibitem{Si2} L. Silvestre, H\"{o}lder estimates for advection fractional-diffusion equations,  arXiv:1009.5723v1 [math.AP] 29 Sep 2010.

\bibitem{Sta} A. Stefanov, Global well-posedness for the 2D quasi-geostrophic equation in a critical Besov space, {\it Electron. J. Differential Equations \bf 2007} (2007),  9 pp.

\bibitem{St} E. Stein, {\it Singular Integrals and Differentiability Properties of Functions}, Princeton Unviersity Press, Princeton, NJ, 1970.


\bibitem{WaJi} H. Wang and H. Jia,  Local well-posedness for the 2D non-dissipative quasi-geostrophic equation in Besov spaces, {\it Nonlinear Anal. \bf 70} (2009), 3791--3798.

\bibitem{WaZh} H. Wang and Z, Zhang,  A frequency localized maximum principle applied to the 2D quasi-geostrophic equation, {\it Comm. Math. Phys.}, in press.


\bibitem{Wu97}J. Wu, Quasi-geostrophic-type equations with initial data in Morrey spaces, {\it Nonlinearity}  {\bf 10} (1997),  1409--1420.

\bibitem{Wu2} J. Wu, Inviscid limits and regularity estimates
for the solutions of the 2-D dissipative quasi-geostrophic
equations, {\it Indiana Univ. Math. J.} {\bf 46} (1997), 1113-1124.

\bibitem{Wu01} J. Wu,  Dissipative quasi-geostrophic equations with $L^p$ data, {\it Electron. J. Differential Equations \bf 2001} (2001), 1-13.

\bibitem{Wu02} J. Wu,  The quasi-geostrophic equation and its two regularizations, {\it Comm. Partial Differential Equations \bf 27} (2002),  1161--1181.

\bibitem{Wu3} J. Wu, Global solutions of the 2D dissipative
quasi-geostrophic equation in
Besov spaces, {\it SIAM J. Math. Anal.}\,\, {\bf 36} (2004/2005),
1014-1030.

\bibitem{Wu4} J. Wu, The quasi-geostrophic equation with critical or supercritical
dissipation, {\it Nonlinearity} \,\,{\bf 18} (2005), 139-154.

\bibitem{Wu41} J. Wu, Solutions of the 2-D quasi-geostrophic equation in H\"{o}lder
spaces, {\it Nonlinear Analysis}\,\, {\bf 62} (2005), 579-594.

\bibitem{Wu31} J. Wu, Lower bounds for an integral involving
fractional Laplacians and the generalized Navier-Stokes equations in
Besov spaces, {\it Comm. Math. Phys.} {\bf 263} (2006), 803-831.

\bibitem{Wu77} J. Wu, Existence and uniqueness results for the 2-D dissipative
quasi-geostrophic equation, {\it Nonlinear Anal.} {\bf 67} (2007),
3013-3036.

\bibitem{Yam} K. Yamazaki, Remarks on the method of modulus of continuity and the modified dissipative Porous Media Equation, {\it J. Differential Equations}, in press.

\bibitem{Yu} X. Yu, Remarks on the global regularity for the super-critical 2D dissipative quasi-geostrophic equation, {\it J. Math. Anal. Appl. \bf 339} (2008), 359--371.

\bibitem{Yuan} B. Yuan, The dissipative quasi-geostrophic equation in weak Morrey spaces, {\it Acta Math. Sin. (Engl. Ser.) \bf 24} (2008), 253--266.

\bibitem{YuanJ} J. Yuan,  On regularity criterion for the dissipative quasi-geostrophic equations, {\it J. Math. Anal. Appl. \bf 340} (2008),  334--339.

\bibitem{Zha0} Z. Zhang,
Well-posedness for the 2D dissipative quasi-geostrophic equations in
the Besov space, {\it Sci. China Ser. A \bf 48} (2005), 1646-1655.

\bibitem{Zha} Z. Zhang,
Global well-posedness for the 2D critical dissipative
quasi-geostrophic equation, {\it Sci. China Ser. A \bf 50} (2007),
485-494.


\bibitem{Zhou} Y. Zhou,
Decay rate of higher order derivatives for solutions to the 2-D
dissipative quasi-geostrophic flows, {\it  Discrete Contin. Dyn.
Syst. \bf 14} (2006), 525-532.

\bibitem{Zhou2} Y. Zhou, Asymptotic behaviour of the solutions to the 2D dissipative quasi-geostrophic flows, {\it Nonlinearity \bf 21} (2008), 2061--2071.

\end{thebibliography}
\end{document}